\documentclass[11pt,a4paper]{article}
\usepackage{amsfonts}
\usepackage{amsmath}
\usepackage{amssymb}
\usepackage{graphicx}
\usepackage{tikz}
\usetikzlibrary{snakes}

\newtheorem{lemma}{\indent Lemma}

\newtheorem{proposition}{\indent Proposition}

\def\kraj{\hfill\rule{6pt}{6pt}}
\def\diag{\mathop{\rm diag}}
\def\deg{\mathop{\rm deg}}
\def\rank{\mathop{\rm rank}}
\def\rk{\mathop{\rm rk}}
\def\sgn{\mathop{\rm sgn}}
\def\lcm{\mathop{\rm lcm}}
\def\tr{\mathop{\rm tr}}
\def\gcd{\mathop{\rm gcd}}
\def\Hom{\mathop{\rm Hom}}
\def\Id{\mathop{\rm Id}}
\def\res{\mathop{\rm res}}
\def\F{\mathbb{F}}
\def\R{\mathbb{R}}
\def\Q{\mathbb{Q}}
\def\Z{\mathbb{Z}}
\def\C{\mathbb{C}}
\def\S{\mathbb{S}}
\def\ka{\mathbb{K}}
\def\E{\mathcal{E}}
\def\U{\mathcal{U}}
\def\B{\mathcal{B}}
\def\Ker{\mathop{\rm Ker}}
\def\max{\mathop{\rm max}}
\def\Seq{\mathbf{\rm Seq}}
\def\Bot{\mathop{\rm Bot}}
\def\Top{\mathop{\rm Top}}
\def\dun{\dot{\U}^+_n}
\def\un{{\U}^+_n}
\def\K{\mathbb{K}}
\def\N{\mathbb{N}}
\def\u3{U^+_q(\mathfrak{sl}_3)}
\def\sl{\mathfrak{sl}}
\date{}
\author{Marko Sto\v si\'c}
\title{On extended graphical calculus for categorified quantum $sl(n)$}

\begin{document}
\maketitle

\begin{abstract}
We study the properties of the extended graphical calculus for categorified quantum $sl(n)$.  
The main results include proofs of Reidemeister 2 and Reidemeister 3-like moves involving strands corresponding to arbitrary thicknesses and arbitrary colors -- the results that were anounced in \cite{canbas}.  
\end{abstract}

\noindent{\small \textbf{Keywords: } categorification; quantum groups; diagrammatic calculus}

\section{Introduction}\label{sec1}

In recent years there has been a lot of work on  diagrammatic categorification of quantum groups, initiated by Lauda's diagrammatic categorification \cite{sl2} (see also  \cite{kl1}) of  Lusztig's idempotented version of $\dot{U}_q(\sl_2)$. This was extended by Khovanov and Lauda in \cite{kl3} to $\dot{U}_q(\sl_n)$ and also in \cite{kl2} to the positive half of an arbitrary quantum group $U^{+}_q(\mathfrak{g})$. 
\\

The general framework of these constructions is to define a certain $2$-category $\U$ whose $1$-morphisms categorify generators of a quantum group, and whose $2$-morphisms are $\ka$-linear combinations of certain planar diagrams modulo local relations, with $\ka$ being a field. Then a $2$-category $\dot{\U}$ is defined as the Karoubi envelope of the $2$-category $\U$, i.e. the smallest category containing $\U$ in which all idempotent $2$-morphisms split. Finally, it is shown that the split Grothendieck group of $\dot{\U}$ is isomorphic to the corresponding quantum group.

For the categorification of ${U}_q^+(\mathfrak{g})$, the  $2$-categories $\U$ and $\dot{\U}$ have a single object. Thus one can see them as monoidal $1$-categories. Since in this paper we are interested in categorifications of positive halves of quantum groups, we shall always assume that $\U$ and $\dot{\U}$ are monoidal ($1$-)categories.\\

The extension of the diagrammatic calculus -- so-called thick calculus -- was introduced in \cite{thick} in the case of quantum $\mathfrak{sl}_2$. With thick calculus one can work directly in $\dot{\U}$, and not just in $\U$. The thick calculus can be extended directly to the case of quantum $\sl_n$ (see \cite{canbas}).

The consequence of \cite{thick}, and of the thick calculus, is that now one can take $\Z$-linear combinations of planar diagrams as morphisms of $\U$, because the idempotents being added have no denominators. 
In this paper, we study the properties of the thick calculus of the category $\dot{\U}^{+}_n$ that categorifies the positive half of the quantum $\sl_n$.  In particular, the main advantage of this approach is that the category $\dot{\U}^{+}_n$ is defined over the ring of integers.\\

Thick calculus already had many applications. This includes the computation of the indecomposable objects of $\dot{\U}^{+}_3$, and consequently the categorification of the Lusztig canonical basis for the positive half of $\dot{U}_q(\sl_3)$ \cite{canbas}. Previous results on this topic were obtained when the category is defined over a field, i.e. when $1$-morphisms are $\ka$-linear combinations of planar diagrams, for some characteristic zero field $\ka$ (see \cite{bk, kl1,vv}).
Furthermore, thick calculus was also used in the combinatorial categorification of $\sl_n$ link invariants of Queffelec and Rose \cite{QR}, via categorical skew-Howe duality.\\

In this paper some of the new relations in the thick calculus of $\sl_n$ are proved. This mainly includes the proofs or Reidemeister 2 and Reidemeister 3 - like moves for strands involving arbitrary thicknesses and labels,  the results that were anounced in \cite{canbas}. Moreover, we also prove some additional relations for passing of a strand through another thick strand labeled by adjacent colors.\\

The results are combinatorial and rely heavily on properties of Schur polynomials, as well as on the numerous properties of the thick diagrammatical calculus both for quantum $\sl_2$ \cite{thick} and quantum $\sl_n$ from \cite{canbas}. In order to make the paper as self-contained as possible, in the first part we recall the basic definitions of the category $\dot{\U}^{+}_n$ and of the thick calculus, together with some basic properties that will be used in the proofs. Mainly we follow \cite{canbas}. Sections \ref{spr2} and \ref{spr3} contain the proofs of the Reidemeister 2 and Reimdeister 3 - like moves, respectively, while in Section \ref{adtr} we prove some additional thick calculus relations as consequences of our main result.\\

\noindent {\bf Acknowledgments}

I am thankful to Mikhail Khovanov on his comments and suggestions on the early version of this paper. This work was partially supported by the Ministry of Education and Science of Serbia, project no. 174012, and  also by the ERC Starting Grant no. 335739 ``Quantum fields and knot homologies" funded by the European Research Council under the European Union's Seventh Framework Programme.
The final part of the work on this project was done while visiting Institute des Hautes Etudes Scientifiques (IHES), at Bures-sur-Yvette, France, and  Max-Planck Institute for Mathematics in Bonn, Germany, and I would like to thank both institutions for hospitality and excellent working atmosphere.

\section{${U}^+_q(\mathfrak{sl}_n)$}\label{sec2}

In this section we define the positive half of quantum 
$\mathfrak{sl}_n$ -- denoted ${U}^+_q(\mathfrak{sl}_n)$. We also give some of its combinatorial properties  in the case $n=3$. 

Let $n\ge 2$ be fixed. The index set of  quantum $\mathfrak{sl}_n$ is $I=\{1,2,\ldots,n-1\}$. 
An inner product is defined on $\Z[I]$ by setting 
\[ 
i\cdot j= \left\{\begin{array}{rl} 2,&\quad i=j\\ -1,&\quad |i-j|=1\\
0,&\quad |i-j|\ge 2\end{array}\right.
\]
for $i,j\in I$.

${U}^+_q(\mathfrak{sl}_n)$ is a  $\Q(q)$-algebra generated by $E_1,E_2,\ldots,E_{n-1}$ modulo relations:
\begin{eqnarray}
E_i^2E_j+E_jE_i^2&=&[2]E_iE_jE_i,\quad i\cdot j=-1,\\
E_iE_j&=&E_jE_i, \quad\quad i\cdot j =0.
\end{eqnarray}

The quantum integers and binomial coefficients are given by:

\begin{eqnarray*}
[n]&=& \frac{q^n-q^{-n}} {q-q^{-1}}, \\
\left[n\right] ! &=& [n][n-1]\cdots[2][1],\\
\left[\!\!\begin{array}{c} n \\ k\end{array}\!\!\right]&=&\frac{[n]!}{[k]![n-k]!}.
\end{eqnarray*}

The divided powers of the generators are defined by 
$$E_i^{(a)}:=\frac{E_i^a}{[a]!}, \quad a\ge 0,\,\,\,\, i=1,\ldots,n-1.$$

The divided powers satisfy:
\begin{eqnarray}
E_i^{(a)}E_j^{(b)}&=&E_j^{(b)}E_i^{(a)},\quad i\cdot j=0,\label{ko}\\
E_{i}^{(a)}E_{i}^{(b)}&=&\left[\begin{array}{c}a+b\\a\end{array}\right]E_{i}^{(a+b)},\label{eiab}
\end{eqnarray}
and  the quantum Serre relations
\begin{equation}
E_i^{(2)}E_j+E_jE_i^{(2)} =  E_iE_jE_i,\quad i\cdot j = -1.\label{kser}
\end{equation}
\vskip 0.3cm
The integral form ${}_{\Z}{U}^+_q(\mathfrak{sl}_n)$ is the 
$\Z[q,q^{-1}]$-subalgebra of ${U}^+_q(\mathfrak{sl}_n)$ generated by $E_{i}^{(a)}$, for all $i=1,\ldots,n-1$ and $a\ge 0$.

\section{The category $\dot{\U}^+_n$}\label{sec3}
A categorification of the positive half of quantum $\mathfrak{sl}_n$ (and also of an arbitrary quantum group 
$U^+_q(\mathfrak{g})$) was defined in \cite{kl2}, though in this paper we prefer the description found in \cite{kl3} in terms of a diagrammatic category  $\mathcal{U}_n^+$.
Before going to the definition of $\U_n^+$, first we recall some notation and explain the diagrams that appear in its definition (see also \cite{kl3}). 

Let $n\ge 2$ be fixed. We refer to the elements of the set $\{1,\ldots,n-1\}$ as \emph{colors}. Let
$\Seq$ denote the set of all finite sequence of colors. For $\nu=(\nu_1,\ldots,\nu_{n-1})\in \N^{n-1}$ let $\Seq(\nu)$ denote the set of all sequences $\underline{i}=(i_1,\ldots,i_k)\in \Seq$ such that $\sharp\{j| i_j=l\}=\nu_l$, for all $l=1,\ldots,n-1$.
Note that $k=\sum_l \nu_l$.\\

We will use the following notion of planar diagrams: We consider 
collections of arcs on the plane connecting the points $\{1,2,\ldots,k\}\subset \R$ in one horizontal line to the points $\{1,2,\ldots,k\}\subset \R$ in another horizontal line. Each arc is labelled by a number from the set $\{1,\ldots,n-1\}$ (called the \emph{color} of the arc). We require that arcs have no critical points when projected to $y$-axis. Arcs can intersect, but no triple intersections are allowed. Finally, an arc can carry dots.

The following is an example of a planar diagram:

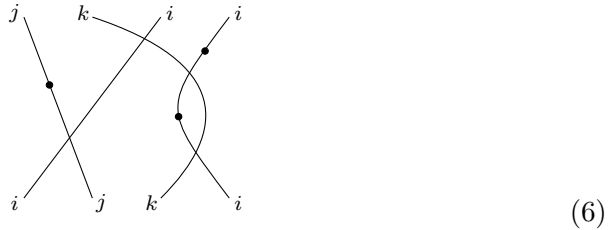
\begin{equation}\label{primer}
\begin{tikzpicture} [scale=0.6]

\draw (0,-2)-- (3,2);
\draw (1.5,-2)--(0,2);

\draw (-0.2,-2.1) node {$\scriptstyle{i}$};
\draw (1.7,-2.1) node {$\scriptstyle{j}$};

\draw (-0.2,2.1) node {$\scriptstyle{j}$};
\draw (3.2,2.1) node {$\scriptstyle{i}$};

\draw (3,-2)..controls (4.5,-0.5) and (4.5,1) ..(1.5,2);
\draw (4.5,-2)..controls (3,0)..(4.5,2);

\draw (2.8, -2.1) node {$\scriptstyle{k}$};
\draw (1.3,2.1) node {$\scriptstyle{k}$};

\draw (4.7,-2.1) node {$\scriptstyle{i}$};
\draw (4.7,2.1) node {$\scriptstyle{i}$};

\filldraw[black] (0.56,0.5) circle (2pt)
(3.39,-0.2) circle (2pt)
(3.97,1.25) circle (2pt);

\end{tikzpicture}
\end{equation}

We identify two planar diagrams if there exists an isotopy between them that does 
not create critical points for the projection onto the $y$-axis.

Since we are not allowing the arcs to have critical points when projected to $y$-axis, we can assume that they are always oriented upwards. We think of a planar diagram as going from its bottom boundary (a sequence of colors) to its top boundary. We read the colors on each boundary from left to right. 

Each diagram has a degree defined as follows.
The degree of a dot is equal to 2. The degree of a crossing between two arcs that are colored  $i$ and $j$ is equal to $-i\cdot j$. In other words, for $i=j$ the degree of a crossing is equal to $-2$, for $|i-j|=1$ (adjacent colors) the degree of a crossing is equal to $1$, while for $|i-j|\ge 2$ (distant colors) the degree of a crossing is equal to $0$. Finally, the degree of a diagram is obtained by summing the contributions coming from all dots and all crossings.

\begin{center}
\begin{tikzpicture} [scale=0.3]
\draw (0,-1.5)-- (0,1.5);

\draw (8,-1.5)--(10,1.5);
\draw (8,1.5)--(10,-1.5);

\draw (4,0) node {$,$};
\draw (-0.2,-1.7) node {$\scriptstyle{i}$};

\draw (7.8,-1.7) node {$\scriptstyle{i}$};
\draw (10.2,-1.7) node {$\scriptstyle{j}$};

\filldraw[black] (0,0) circle (4pt);

\draw (-6,-3.2) node {degree:};

\draw (0,-3.2)  node {$\scriptstyle{+2}$};

\draw (8.8,-3.2)  node {$\scriptstyle{-i\cdot j}$};

\end{tikzpicture}
\end{center}

We also use the following shorthand for a collection of  dots on a strand.
\begin{center}
\begin{tikzpicture} [scale=0.35]

\draw (0,-2)-- (0,2);
\draw (-0.2,-2.1) node {$\scriptstyle{i}$};
\draw (-0.4,0) node {$\scriptstyle{d}$};

\draw (1.2,0) node {$ := $};

\draw (3,-2) --(3,2);

\draw (2.8,-2.1) node {$\scriptstyle{i}$};
\draw (2.73,0.55) node {$\vdots$};

\draw (4.5,0) node {$ d $};

\draw  [thick] (3.5,1.1)..controls(3.7,1.1)..(3.7,0.2);
\draw  [thick] (3.7,0.2)--(3.9,0);
\draw  [thick] (3.7,-0.2)--(3.9,0);
\draw  [thick] (3.5,-1.1)..controls(3.7,-1.1)..(3.7,-0.2);

\filldraw[black] (0,0) circle (2pt)
(3,-1) circle (2pt)
(3,1) circle (2pt)
(3,-0.5) circle (2pt);

\end{tikzpicture}
\end{center}

\subsection{The category $\U^+_n$}

$\U^+_n$ is the monoidal $\Z$-linear additive category whose objects and morphisms are the following:\\

$\bullet$    objects: for each $\underline{i}=(i_1,\ldots,i_k)\in \Seq$ and $t\in\Z$, we define $\E_{\underline{i}}\{t\}:=\E_{i_1}\ldots \E_{i_k}\{t\}$. An object of $\U^+_n$ is a formal finite direct sum of $\E_{\underline{i}}\{t\}$, with $i\in \Seq$ and $t\in\Z$.\\

$\bullet$    morphisms: for $\underline{i}=(i_1,\ldots,i_k)\in \Seq({\nu})$ and $\underline{j}=(j_1,\ldots,j_l)\in \Seq(\mu)$ the 
set $\Hom(\E_{\underline{i}}\{t\},\E_{\underline{j}}\{t'\})$ is empty, unless $\nu=\mu$. If $\nu=\mu$ (and consequently $k=l$), the 
morphisms from $\E_{\underline{i}}\{t\}$ to $\E_{\underline{j}}\{t'\}$ consist of finite $\Z$-linear combinations of planar diagrams going from $\underline{i}$ to $\underline{j}$, of degree $t-t'$, modulo the following set of homogeneous local relations:







\begin{center}
\begin{tikzpicture} [scale=0.3]

\draw (6,-2)-- (8,2);
\draw (8,-2)--(6,2);

\draw (5.7,-2) node {$\scriptstyle{i}$};
\draw (8.3,-2) node {$\scriptstyle{i}$};

\draw (9,0) node {$\displaystyle{-}$};

\filldraw[black] (6.5,1) circle (3pt);

\draw (10,-2)-- (12,2);
\draw (12,-2)--(10,2);

\draw (9.7,-2) node {$\scriptstyle{i}$};
\draw (12.3,-2) node {$\scriptstyle{i}$};

\filldraw[black] (11.5,-1) circle (3pt);

\draw (13,0) node {$\displaystyle{=}$};

\draw (14,-2)-- (16,2);
\draw (16,-2)--(14,2);

\draw (13.7,-2) node {$\scriptstyle{i}$};
\draw (16.3,-2) node {$\scriptstyle{i}$};

\filldraw[black] (14.5,-1) circle (3pt);

\draw (17,0) node {$\displaystyle{-}$};

\draw (18,-2)-- (20,2);
\draw (20,-2)--(18,2);

\draw (17.7,-2) node {$\scriptstyle{i}$};
\draw (20.3,-2) node {$\scriptstyle{i}$};

\filldraw[black] (19.5,1) circle (3pt);

\draw (21,0) node {$\displaystyle{=}$};

\draw (22,-2)-- (22,2);
\draw (23.5,-2)--(23.5,2);

\draw (21.8,-2) node {$\scriptstyle{i}$};
\draw (23.7,-2) node {$\scriptstyle{i}$};

\end{tikzpicture}
\end{center}

\begin{center}
\begin{tikzpicture} [scale=0.3]

\draw (-8,-2)..controls (-6,0) ..(-8,2);
\draw (-6,-2)..controls (-8,0)..(-6,2);
\draw (-8.3,-2) node {$\scriptstyle{i}$};
\draw (-5.7,-2) node {$\scriptstyle{i}$};

\draw (-5,0) node {$\displaystyle{= 0} , $};

\draw (3.5,-2)--(6.5,2);
\draw (3.5,2)--(6.5,-2);
\draw (5,-2)..controls (3.5,0) ..(5,2);

\draw (3.7,-2) node {$\scriptstyle{i}$};
\draw (5.3,-2) node {$\scriptstyle{i}$};
\draw (6.8,-2) node {$\scriptstyle{i}$};

\draw (7.5,0) node {$=$};

\draw (8.5,-2)--(11.5,2);
\draw (8.5,2)--(11.5,-2);
\draw (10,-2)..controls (11.5,0) ..(10,2);

\draw (8.7,-2) node {$\scriptstyle{i}$};
\draw (10.3,-2) node {$\scriptstyle{i}$};
\draw (11.8,-2) node {$\scriptstyle{i}$};

\end{tikzpicture}
\end{center}

\begin{equation}\label{r2tan}
\begin{tikzpicture} [scale=0.3]
\draw (-1,-2)..controls (1,0) ..(-1,2);
\draw (1,-2)..controls (-1,0)..(1,2);
\draw (-1.3,-2) node {$\scriptstyle{i}$};
\draw (1.3,-2) node {$\scriptstyle{j}$};

\draw (2,0) node {$\displaystyle{=}$};

\draw (3,-2)--(3,2);
\draw (4.5,-2)--(4.5,2);

\draw (5.5,0) node {$+$};

\draw (6.5,-2)--(6.5,2);
\draw (8,-2)--(8,2);

\filldraw[black] (3,0) circle (3pt)
(8,0) circle (3pt);

\draw (15,0) node {,$\quad$ when $\quad i\cdot j=-1$};

\end{tikzpicture}
\end{equation}

\begin{center}
\begin{tikzpicture}[scale=0.3]

\draw (11,-2)..controls (13,0) ..(11,2);
\draw (13,-2)..controls (11,0)..(13,2);
\draw (10.7,-2) node {$\scriptstyle{i}$};
\draw (13.3,-2) node {$\scriptstyle{j}$};

\draw (14,0) node {$\displaystyle{=}$};

\draw (15,-2)--(15,2);
\draw (16.5,-2)--(16.5,2);

\draw (14.7,-2) node {$\scriptstyle{i}$};
\draw (16.8,-2) node {$\scriptstyle{j}$};

\draw (23,0) node {,$\quad$ when $\quad i\cdot j=0$};

\end{tikzpicture}
\end{center}

\begin{equation}\label{thinslide}
\begin{tikzpicture} [scale=0.3]
\draw (6,-2)-- (8,2);
\draw (8,-2)--(6,2);

\draw (5.7,-2) node {$\scriptstyle{i}$};
\draw (8.3,-2) node {$\scriptstyle{j}$};

\draw (9,0) node {$\displaystyle{=}$};

\filldraw[black] (6.5,1) circle (3pt);

\draw (10,-2)-- (12,2);
\draw (12,-2)--(10,2);

\draw (9.7,-2) node {$\scriptstyle{i}$};
\draw (12.3,-2) node {$\scriptstyle{j}$};

\filldraw[black] (11.5,-1) circle (3pt);

\draw (15,0) node {and};


\draw (19,-2)-- (21,2);
\draw (21,-2)--(19,2);

\draw (18.7,-2) node {$\scriptstyle{i}$};
\draw (21.3,-2) node {$\scriptstyle{j}$};

\filldraw[black] (19.5,-1) circle (3pt);

\draw (22,0) node {$\displaystyle{=}$};

\draw (23,-2)-- (25,2);
\draw (25,-2)--(23,2);

\draw (22.7,-2) node {$\scriptstyle{i}$};
\draw (25.3,-2) node {$\scriptstyle{j}$};

\filldraw[black] (24.5,1) circle (3pt);

\draw (31.5,0) node {, $\quad$ when $\quad i\ne j$};

\end{tikzpicture}
\end{equation}

\begin{center}
\begin{tikzpicture} [scale=0.4]

\draw (-1.5,-2)--(1.5,2);
\draw (-1.5,2)--(1.5,-2);
\draw (0,-2)..controls (-1.5,0) ..(0,2);

\draw (-1.3,-2) node {$\scriptstyle{i}$};
\draw (0.3,-2) node {$\scriptstyle{j}$};
\draw (1.8,-2) node {$\scriptstyle{k}$};

\draw (2.5,0) node {$=$};

\draw (3.5,-2)--(6.5,2);
\draw (3.5,2)--(6.5,-2);
\draw (5,-2)..controls (6.5,0) ..(5,2);

\draw (3.7,-2) node {$\scriptstyle{i}$};
\draw (5.3,-2) node {$\scriptstyle{j}$};
\draw (6.8,-2) node {$\scriptstyle{k}$};

\draw (14,0) node {, if $i\ne k$ or $i\cdot j\ne -1$};

\end{tikzpicture}
\end{center}

\begin{equation}\label{r3tan}
\begin{tikzpicture} [scale=0.4]

\draw (-1.5,-2)--(1.5,2);
\draw (-1.5,2)--(1.5,-2);
\draw (0,-2)..controls (-1.5,0) ..(0,2);

\draw (-1.3,-2) node {$\scriptstyle{i}$};
\draw (0.3,-2) node {$\scriptstyle{j}$};
\draw (1.8,-2) node {$\scriptstyle{i}$};

\draw (2.5,0) node {$=$};

\draw (3.5,-2)--(6.5,2);
\draw (3.5,2)--(6.5,-2);
\draw (5,-2)..controls (6.5,0) ..(5,2);

\draw (3.7,-2) node {$\scriptstyle{i}$};
\draw (5.3,-2) node {$\scriptstyle{j}$};
\draw (6.8,-2) node {$\scriptstyle{i}$};

\draw (8,0) node {$+$};

\draw (9,-2)--(9,2);
\draw (10.5,-2)--(10.5,2);
\draw (12,-2)--(12,2);

\draw (9.3,-2) node {$\scriptstyle{i}$};
\draw (10.8,-2) node {$\scriptstyle{j}$};
\draw (12.3,-2) node {$\scriptstyle{i}$};

\draw (18,0) node {, if $i\cdot j= -1$};

\end{tikzpicture}
\end{equation}

This ends the definition of $\U_n^+$. \\

As an example of a morphism, a diagram from (\ref{primer}) represents a morphism in $\Hom(\E_i\E_j\E_k\E_i\{t\},\E_j\E_k\E_i\E_i\{t+2\})$.\\

We have the following relation in $\U_n^+$:

\begin{proposition}[Dot Migration]\label{dotm}\cite[Proposition 5.2]{sl2}
We have

\begin{center}
\begin{tikzpicture} [scale=0.4]

\draw (6,-2)-- (8,2);
\draw (8,-2)--(6,2);

\draw (5.7,-2) node {$\scriptstyle{i}$};
\draw (8.3,-2) node {$\scriptstyle{i}$};

\draw (9,0) node {$\displaystyle{-}$};

\filldraw[black] (6.5,1) circle (3pt);
\draw (6.1,1) node {$\scriptstyle{d}$};

\draw (10,-2)-- (12,2);
\draw (12,-2)--(10,2);

\draw (9.7,-2) node {$\scriptstyle{i}$};
\draw (12.3,-2) node {$\scriptstyle{i}$};

\filldraw[black] (11.5,-1) circle (3pt);
\draw (11.9,-1) node {$\scriptstyle{d}$};

\draw (13,0) node {$\displaystyle{=}$};

\draw (14,-2)-- (16,2);
\draw (16,-2)--(14,2);

\draw (13.7,-2) node {$\scriptstyle{i}$};
\draw (16.3,-2) node {$\scriptstyle{i}$};

\filldraw[black] (14.5,-1) circle (3pt);
\draw (14.1,-1) node {$\scriptstyle{d}$};

\draw (17,0) node {$\displaystyle{-}$};

\draw (18,-2)-- (20,2);
\draw (20,-2)--(18,2);

\draw (17.7,-2) node {$\scriptstyle{i}$};
\draw (20.3,-2) node {$\scriptstyle{i}$};

\filldraw[black] (19.5,1) circle (3pt);
\draw (19.9,1) node {$\scriptstyle{d}$};

\draw (22.5,0) node {$\displaystyle{=} \sum_{r+s=d-1}$};

\draw (25.5,-2)-- (25.5,2);
\draw (27,-2)--(27,2);

\draw (25.3,-2) node {$\scriptstyle{i}$};
\draw (27.2,-2) node {$\scriptstyle{i}$};

\filldraw[black] (25.5,0) circle (3pt);
\draw (25.2,0) node {$\scriptstyle{r}$};

\filldraw[black] (27,0) circle (3pt);
\draw (27.4,0) node {$\scriptstyle{s}$};

\end{tikzpicture}
\end{center}

\end{proposition}

\subsection{The category $\dot{\U}_n^+$ and thick calculus}

In \cite{thick}, the extension of the calculus to thick edges have been introduced. Thick lines categorify the divided powers $E_i^{(a)}$ (see below and Section 4 of \cite{thick}). \\

For a category $\mathcal{C}$, the Karoubi envelope $Kar(\mathcal{C})$ is the smallest category containing $\mathcal{C}$, such that all idempotents split (for more details, see e.g. Section 3.4 of \cite{thick}).

We define the category $\dot{\U}^+_n$ as the Karoubi envelope of the category $\U^+_n$. \\

As in \cite{thick}, the category $\dot{\U}^+_n$ categorifies $U_q^+(\mathfrak{sl}_n)$, in a sense that its split Grothendieck group is isomorphic to the integral form of  $U_q^+(\mathfrak{sl}_n)$. The isomorphism sends the class of  $\E_i^{(a)}$ to the generator $E_i^{(a)}$ of $U_q^+(\mathfrak{sl}_n)$.\\

In the category $\dot{\U}^+_n$, the planar diagrams with thick edges from above can be interpreted as morphisms whose bottom and top end correspond to certain objects of $\dot{\U}^+_n$. In particular, the object corresponding to bottom (or the top end)  of an arc of color $i$ and thickness $a$ is denoted $\E_i^{(a)}$. 

A thick line of color $i$ is defined as the identity morphism:

\begin{center}
\begin{tikzpicture} [scale=0.3]
\draw [thick]  (0,-1.5)-- (0,1.5);

\draw (0.5,-1.7) node {$\scriptstyle{a}$};

\draw (1.5,-0.1) node {$\scriptstyle{:}$};

\draw (6,0.1) node {$\E_i^{(a)}\longrightarrow \E_i^{(a)}$};

\end{tikzpicture}
\end{center}

It is given explicitly in terms of ``ordinary" lines from above - see \cite[Equation (2.18)]{thick}, and drawn as a strand with an additional label (natural number) $a$ (also called the \emph{thickness} of a strand). In particular, the ordinary strands from above correspond to the case $a=1$, and are also called \emph{thin} edges or thin strands. 
We refer the reader to \cite{thick}, in particular Sections 2 and 4, for more details. Here we just recall the basic facts that will be used later on.\\

Trivalent vertices of a single color are now allowed in our planar diagrams, as long as the sum of thicknesses of the incoming edges is equal to the sum of  thicknesses of the outgoing edges. These are called splitters in \cite{thick}. The trivalent vertices (for any color $i$ - the labels on the pictures below represent thicknesses)

\begin{center}
\begin{tikzpicture} [scale=0.4]

\draw[thick] (0,-1)-- (1,0);
\draw[thick] (2,-1)--(1,0);
\draw[thick] (1,0)--(1,1);

\draw (-0.25,-1) node {$\scriptstyle{a}$};
\draw (2.25,-1) node {$\scriptstyle{b}$};
\draw (0,1) node {$\scriptstyle{a+b}$};

\draw (7,0) node {$:\,\, \E_i^{(a)}\E_i^{(b)}\rightarrow\E_i^{(a+b)}$};

\draw[thick] (15,-1)-- (15,0);
\draw[thick] (14,1)--(15,0);
\draw[thick] (15,0)--(16,1);

\draw (15.8,-1) node {$\scriptstyle{a+b}$};
\draw (13.65,1) node {$\scriptstyle{a}$};
\draw (16.25,1) node {$\scriptstyle{b}$};

\draw (21,0) node {$:\,\,\E_i^{(a+b)} \rightarrow\E_i^{(a)}\E_i^{(b)}$};

\end{tikzpicture}
\end{center}






\noindent correspond to the projection and inclusion maps, respectively, obtained from the decomposition:
\[
\E_i^{(a)}\E_i^{(b)}\cong \bigoplus_{\left[a+b\atop b\right]}\E_i^{(a+b)}
\]

The degrees of both of this two vertices are equal to $-ab$, which explains which summands must be involved in these morphisms. The explicit definitions are given in \cite[pp. 15]{thick}.

These morphisms may be composed, and in particular they can be used to define the thick crossing:

\begin{center}
\begin{tikzpicture}[scale=0.5]

\draw [thick](5,-2.2)--(6,-1.3);
\draw [thick](6,-1.3)--(6,1.3);
\draw [thick](6,1.3)--(5,2.2);
\draw [thick](7,-2.2)--(6,-1.3);
\draw [thick](6,1.3)--(7,2.2);
\draw (3.5,0) node {$=$};

\draw [thick](0,-2.2)--(2,2.2);
\draw [thick](0,2.2)--(2,-2.2);


\draw (4.65,2.3) node {$\scriptstyle{b}$};
\draw (4.65,-2.3) node {$\scriptstyle{a}$};
\draw (5.3,0) node {$\scriptstyle{a+b}$};
\draw (7.2,2.3) node {$\scriptstyle{a}$};
\draw (7.2,-2.3) node {$\scriptstyle{b}$};

\draw (-0.3,2.3) node {$\scriptstyle{b}$};
\draw (-0.3,-2.3) node {$\scriptstyle{a}$};
\draw (2.2,2.3) node {$\scriptstyle{a}$};
\draw (2.2,-2.3) node {$\scriptstyle{b}$};

\end{tikzpicture}
\end{center}


\vskip 0.5cm

\subsubsection{Some properties of the thick calculus}

Below we give some of the basic properties of thick edges that we shall use in this paper (see \cite{thick} for more details).
Note that the labels of the strands below denote thickness. All relations hold under
horizontal and vertical flips, because of the symmetries on $\U$.

\begin{proposition}[Associativity of splitters]\label{asos}\cite[Proposition 2.2.4]{thick}
For arbitrary color $i$  we have the following:

\begin{center}
\begin{tikzpicture}[scale=0.4]
\draw [thick] (-2,2)--(0,-0.5);
\draw [thick](-1,0.75)--(0,2);
\draw [thick](0,-0.5)--(2,2);
\draw [thick](0,-0.5)--(0,-2);

\draw (-2.3,2) node {$\scriptstyle{a}$};
\draw (-0.45,2) node {$\scriptstyle{b}$};
\draw (2.3,2) node {$\scriptstyle{c}$};
\draw (-1.3,0) node {$\scriptstyle{a+b}$};
\draw (1.1,-2.2) node {$\scriptstyle{a+b+c}$};

\draw (3,0) node {$=$};

\draw [thick](4,2)--(6,-0.5);
\draw [thick](7,0.75)--(6,2);
\draw [thick](6,-0.5)--(8,2);
\draw [thick](6,-0.5)--(6,-2);

\draw (3.7,2) node {$\scriptstyle{a}$};
\draw (5.7,2) node {$\scriptstyle{b}$};
\draw (8.3,2) node {$\scriptstyle{c}$};
\draw (7.4,0) node {$\scriptstyle{b+c}$};
\draw (7.5,-2.2) node {$\scriptstyle{a+b+c}$};

\end{tikzpicture}
\end{center}

\end{proposition}

In particular, the above Associativity of Splitters imply that the multi-splitters,
i.e. splitting of a thick line into three or more strands, is well-defined.

\begin{center}
\begin{tikzpicture}[scale=0.4]
\draw [thick] (-2,2)--(0,-0.5);
\draw [thick](0,-0.5)--(0,2);
\draw [thick](0,-0.5)--(2,2);
\draw [very thick](0,-0.5)--(0,-2);

\draw (-2.3,2) node {$\scriptstyle{a}$};
\draw (-0.45,2) node {$\scriptstyle{b}$};
\draw (2.3,2) node {$\scriptstyle{c}$};
\draw (1.1,-2.2) node {$\scriptstyle{a+b+c}$};

\draw (3,0) node {$=$};

\draw [thick](4,2)--(6,-0.5);
\draw [thick](7,0.75)--(6,2);
\draw [thick](6,-0.5)--(8,2);
\draw [very thick](6,-0.5)--(6,-2);

\draw (3.7,2) node {$\scriptstyle{a}$};
\draw (5.7,2) node {$\scriptstyle{b}$};
\draw (8.3,2) node {$\scriptstyle{c}$};
\draw (7.4,0) node {$\scriptstyle{b+c}$};
\draw (7.5,-2.2) node {$\scriptstyle{a+b+c}$};

\end{tikzpicture}
\end{center}

\begin{proposition}[Pitchfork lemma]\cite{thick}
For any two colors $i$ (drawn as a thick line) and $j$ (drawn dashed) we have:
\begin{center}
\begin{tikzpicture}[scale=0.5]
\draw [thick](-2,-2)--(-1,-1);
\draw [thick](-1,-1)..controls (-0.5,0.5)..(1,2);
\draw [thick](-1,-1)..controls (0.5,-0.5)..(2,1);
\draw [thick,dashed] (-2,2)--(2,-2);

\draw (1.3,2) node {$\scriptstyle a$};
\draw (-2,-2.35) node {$\scriptstyle a+b$};
\draw (2.4,1.2) node {$\scriptstyle b$};
\draw (2.3,-2) node {$\scriptstyle c$};

\draw (3,0) node {$=$};

\draw [thick](4,-2)--(6.5,0.5);
\draw [thick](6.5,0.5)..controls (6.6,1.5)..(7,2);
\draw [thick](6.5,0.5)..controls (7.5,0.6)..(8,1);

\draw [thick,dashed] (4,2)--(8,-2);

\draw (7.4,2) node {$\scriptstyle a$};
\draw (4.5,-2.3) node {$\scriptstyle a+b$};
\draw (8.35,1.2) node {$\scriptstyle b$};
\draw (8.3,-2) node {$\scriptstyle c$};

\end{tikzpicture}
\end{center}

\end{proposition}

\begin{proposition}[Opening of a Thick Edge]\cite[Proposition 2.2.5]{thick}\label{OTE}
For any color $i$ we have:

\begin{center}
\begin{tikzpicture}[scale=0.9]

\draw [thick](0,-1.5)--(1,-0.6);
\draw [thick](1,-0.6)--(1,0.6);
\draw [thick](1,0.6)--(0,1.5);
\draw [thick](2,-1.5)--(1,-0.6);
\draw [thick](1,0.6)--(2,1.5);
\draw (3,0) node {$=$};

\draw [thick](4,-1.5)--(7,1.5);
\draw [thick](4,1.5)--(7,-1.5);
\draw [thick](4.6,-0.91)--(4.6,0.91);

\draw (-0.5,1.6) node {$\scriptstyle{b+x}$};
\draw (-0.5,-1.6) node {$\scriptstyle{a+x}$};
\draw (0.3,0) node {$\scriptstyle{a+b+x}$};
\draw (2.2,1.6) node {$\scriptstyle{a}$};
\draw (2.2,-1.6) node {$\scriptstyle{b}$};

\draw (3.5,1.6) node {$\scriptstyle{b+x}$};
\draw (3.5,-1.6) node {$\scriptstyle{a+x}$};
\draw (4.3,0) node {$\scriptstyle{x}$};
\draw (7.2,1.6) node {$\scriptstyle{a}$};
\draw (7.2,-1.55) node {$\scriptstyle{b}$};
\draw (5,0.7) node {$\scriptstyle{b}$};
\draw (5,-0.7) node {$\scriptstyle{a}$};

\end{tikzpicture}
\end{center}

\end{proposition}

\subsection{Schur polynomials and decorations of thick lines}

Just as we can decorate thin strands with dots, we can decorate thick lines with symmetric polynomials. These correspond to symmetric polynomials in dots on thin edges involved in the definition of a thick line (for a precise definition, see \cite{thick}). 
For notational convenience we will only decorate thick strands with Schur polynomials, which form an additive basis of the ring of symmetric polynomials. 

\subsubsection{Schur polynomials}
Here we recall briefly the definition and some basic notation and properties of  Schur polynomials. For more details, see e.g. \cite{fulton,thick,mac}.

By a partition $\alpha=(\alpha_1,\ldots,\alpha_k)$, we mean a non-increasing sequence of non-negative integers. We identify two partitions if they differ by a sequence of zeros at the end. We set $|\alpha|=\sum_i \alpha_i$. If for some $a$ we have $\alpha_{a+1}=0$, we say that $\alpha$ has at most $a$ parts. We denote the set of all partitions with at most $a$ parts by 
$P(a)$. Furthermore, by $P(a,b)$ we denote the subset of all partitions $\alpha$ from $P(a)$ such that $\alpha_1\le b$. In other words, $P(a,b)$ consists of partitions fitting inside a rectangle with $a$ rows and $b$ columns. The partition corresponding to this rectangle we denote by $K_{a,b}$, i.e. $K_{a,b}=(\underbrace{b,b,\ldots,b}_a)$.

We shall need to express quantum binomial coefficients as a sum over partitions fitting inside a rectangle. For any two nonnegative integers $a$ and $b$ we have:

\begin{equation}\label{qupa}
\left[\begin{array}{c} a+b \\ a \end{array}\right]=\sum_{\alpha\in P(a,b)} q^{2|\alpha|-ab}.
\end{equation}

By $\bar{\alpha}$ we denote the dual (conjugate) partition of $\alpha$, i.e. $\alpha_j=\sharp\{i| \alpha_i\ge j\}$. If $\alpha\in P(a,b)$, we define partition $\hat{\alpha}$ by $\hat{\alpha}=\overline{(b-\alpha_a,\ldots,b-\alpha_1)}$. Note that if $\alpha\in P(a,b)$, then $\bar{\alpha}\in P(b,a)$ and $\hat{\alpha}\in P(b,a)$.

For any partition $\alpha\in P(a)$, the Schur polynomial $\pi_{\alpha}$ is given by the formula:
\[
\pi_{\alpha}(x_1,x_2,\ldots,x_a)=\frac{|x_i^{\alpha_j+a-j}|}{\Delta},
\]
where $\Delta=\prod_{1\le r \le s\le a} (x_r-x_s)$, and $|x_i^{\alpha_j+a-j}|$ is the determinant of the $a\times a$ matrix whose $(i,j)$ entry is $x_i^{\alpha_j+a-j}$. We extend our notation, so that $\pi_{\alpha}(x_1,x_2,\ldots,x_a)=0$ is some entry of $\alpha$ is negative ($\alpha$ is not a partition then), or if $\alpha_{a+1}>0$.

For two partitions $\alpha$ and $\gamma$, we say that $\alpha\subset\gamma$ if $\alpha_i\le \gamma_i$ for all $i\ge 1$. 

For three partitions $\alpha$, $\beta$ and $\gamma$, the Littlewood-Richardson coefficients $c_{\alpha,\beta}^{\gamma}$ are given by:
\[
\pi_{\alpha}\pi_{\beta}=\sum_{\gamma} c_{\alpha,\beta}^{\gamma} \pi_{\gamma}.
\]
The coefficients $c_{\alpha,\beta}^{\gamma}$ are nonnegative integers that can be nonzero only when $|\gamma|=|\alpha|+|\beta|$. Also,  $c_{\alpha,\beta}^{\gamma}\ne 0$ only when $\alpha \subset \gamma$ and $\beta \subset \gamma$. In particular:
\begin{equation}\label{LRpr}
\alpha\in P(a,x), \beta\in P(b,y)\textrm{ and }c_{\alpha,\beta}^{\gamma}\ne 0,\textrm{ imply } \gamma\in P(a+b,x+y).
\end{equation}

The Littlewood-Richardson coefficients can be naturally extended for more than three partitions: for partitions $\alpha_1,\ldots,\alpha_k$ and $\beta$, with $k\ge 2$, we define $c_{\alpha_1,\ldots,\alpha_k}^{\beta}$ by:
\[
\pi_{\alpha_1}\ldots\pi_{\alpha_k}=\sum_{\beta} c_{\alpha_1,\ldots,\alpha_k}^{\beta} \pi_{\beta}.
\]

For two partitions $\alpha$ and $\gamma$, we define skew-Schur polynomial $\pi_{^{\gamma}/_{\alpha}}$ by:
\[
\pi_{^{\gamma}/_{\alpha}}=\sum_{\beta} c_{\alpha,\beta}^{\gamma} \pi_{\beta}.
\]
It can be nonzero only when $\alpha \subset \gamma$.

If $\gamma=(\gamma_1,\ldots,\gamma_a) \subset K_{a,b}$, then by $K_{a,b}-\gamma$ we denote the partition $(b-\gamma_a,\ldots,b-\gamma_1)$.
For a partition $\nu\in P(a)$, by $\nu+K_{a,b}$ we denote the partition $(\nu_1+b,\ldots,\nu_a+b)$. Furthermore, for every two partitions $\psi\in P(a)$ and $\gamma\in P(a,b)$, we have that $c_{\gamma,\psi}^{\nu+K_{a,b}} = c_{\nu, K_{a,b}-\gamma}^{\psi}$. In particular, if $\nu=\emptyset$, one has $c_{\gamma,\psi}^{K_{a,b}} = c_{\emptyset, K_{a,b}-\gamma}^{\psi}=\delta_{\psi,K_{a,b}-\gamma}$, and so $\pi_{{}_{K_{a,b}/_{\gamma}}}=\pi_{{}_{K_{a,b}-\gamma}}$.

The elementary symmetric polynomials $\varepsilon_m(x_1,\ldots,x_a)$, for $m=0,\ldots,a$ are special Schur polynomials: $\varepsilon_m(x_1,\ldots,x_a)=\pi_{{\scriptsize{(\underbrace{1,1,\ldots,1}_m)}}} (x_1,\ldots,x_a)$. For $m<0$ or $m>a$, we have $\varepsilon_m(x_1,\ldots,x_a)=0$. 

The Schur polynomials can be conveniently expressed as a determinant of a matrix whose entries are the elementary symmetric  polynomials, by the following Giambelli formula: for a partition $\alpha=(\alpha_1,\ldots,\alpha_a)$, we have
\begin{equation}\label{djam}
\pi_{\bar{\alpha}}=\det [\varepsilon_{\alpha_i+j-i}]_{i,j=1}^a.
\end{equation}

\subsubsection{Decorated thick edges}
Here we recall some of the basic properties of decorated thick lines that we shall need in this paper. For more details see \cite{thick}.

A thick line of thickness $a$ can be decorated with any Schur polynomial $\pi_{\alpha}$. For $\alpha\notin P(a)$, the resulting morphism is zero. For $\alpha\in P(a)$ of the form 
$\alpha=(\alpha_1,\alpha_2,\ldots,\alpha_a)$, one can express the decoration of a thick line in terms of thin lines and dots as follows
\begin{equation}\label{decdef}
\begin{tikzpicture}[scale=0.7]

\draw [very thick] (0,1.5)--(0,-1.5);

\draw (-0.5,0) node {$\scriptstyle{\pi_{\alpha}}$};
\draw (2.1,0.15) node {$\scriptscriptstyle{\alpha_{1}+a-1}$};
\draw (4.25,0.2) node {$\scriptscriptstyle{\alpha_{2}\!+\!a\!-\!2}$};
\draw (5.65,-0.3) node {$\scriptscriptstyle{\alpha_{a-1}\!+\!1}$};
\draw (7.7,0) node {$\scriptstyle{\alpha_{a}}$};
\draw (5,0) node {$\cdots$};


\filldraw[black] (0,0) circle (1.8pt)
                  
                 (6.5,0) circle (1.8pt)
                 (7.25,0) circle (1.8pt)
                 (3.522,0) circle (1.8pt)
                 (2.725,0) circle (1.8pt);

\draw (0.2,-1.5) node {$\scriptstyle{a}$};

\draw (5.2,2) node {$\scriptstyle{a}$};
\draw (5.2,-2) node {$\scriptstyle{a}$};

\draw (1,0) node{$=$};
\draw [very thick] (5,-2)--(5,-1.5);
\draw [very thick] (5,2)--(5,1.5);
\draw (5,1.5)..controls (2,0.75) and (2,-0.75)..(5,-1.5);
\draw (5,1.5)..controls (8,0.75) and (8,-0.75)..(5,-1.5);
\draw (5,1.5)..controls (3,0.75) and (3,-0.75)..(5,-1.5);
\draw (5,1.5)..controls (7,0.75) and (7,-0.75)..(5,-1.5);
\end{tikzpicture}
.
\end{equation}

From the definition of the Littlewood-Richardson coefficient, we have

\begin{center}
\begin{tikzpicture}[scale=0.3]

\draw [very thick] (0,-2)--(0,2);

\filldraw[black] (0,-0.8) circle (4pt)
                 (0,0.8) circle (4pt)
                 (9,0) circle (4pt);

\draw (-0.9,-0.8) node {$\scriptstyle{\pi_{\alpha}}$};
\draw (-0.9,0.8) node {$\scriptstyle{\pi_{\beta}}$};
\draw (8.15,0) node {$\scriptstyle{\pi_{\gamma}}$};

\draw (4,0) node {${=\,\sum\nolimits_{\gamma}\,\,c_{\alpha\beta}^{\gamma}}$};

\draw [very thick] (9,-2)--(9,2);

\end{tikzpicture}
\end{center}

By ``exploding" a thick edge into thin edges, we obtain diagrams that are antisymmetric with respect to the exchange of dots on two neighboring strands:
\begin{center}
\begin{tikzpicture}[scale=0.5]

\draw [very thick] (11,-2)--(11,-1.2);
\draw [very thick] (11,2)--(11,1.2);
\draw (11,1.2)..controls (10,0)..(11,-1.2);
\draw (11,1.2)..controls (12,0)..(11,-1.2);
\draw (10,0) node {$\scriptstyle{a}$};
\draw (12,0) node {$\scriptstyle{b}$};
\draw (13,0) node {$= - $};
\draw [very thick] (15,-2)--(15,-1.2);
\draw [very thick] (15,2)--(15,1.2);
\draw (15,1.2)..controls (14,0)..(15,-1.2);
\draw (15,1.2)..controls (16,0)..(15,-1.2);
\draw (14,0) node {$\scriptstyle{b}$};
\draw (16,0) node {$\scriptstyle{a}$};

\filldraw[black] (10.23,0) circle (1.7pt)
                 (11.75,0) circle (1.7pt)
                 (14.23,0) circle (1.7pt)
                 (15.75,0) circle (1.7pt);

\end{tikzpicture}
\end{center}

This antisymmetry implies the following:

\begin{lemma} {}$\quad$

\begin{center}

\begin{tikzpicture}[scale=0.7]
\draw (2.15,0) node {$\scriptscriptstyle{x_1}$};
\draw (3.95,0) node {$\scriptscriptstyle{x_2}$};
\draw (5.85,0) node {$\scriptscriptstyle{x_{a\!-\!1}}$};
\draw (7.6,0) node {$\scriptstyle{x_{a}}$};
\draw (4.96,0) node {$\cdots$};


\filldraw[black] 
                  (7.23,0) circle (1.8pt)
                 (6.5,0) circle (1.8pt)
                 
(3.522,0) circle (1.8pt)
                 (2.725,0) circle (1.8pt);


\draw (5.2,2) node {$\scriptstyle{a}$};
\draw (5.2,-2) node {$\scriptstyle{a}$};

\draw [very thick] (5,-2)--(5,-1.5);
\draw [very thick] (5,2)--(5,1.5);
\draw (5,1.5)..controls (2,0.75) and (2,-0.75)..(5,-1.5);
\draw (5,1.5)..controls (8,0.75) and (8,-0.75)..(5,-1.5);
\draw (5,1.5)..controls (3,0.75) and (3,-0.75)..(5,-1.5);
\draw (5,1.5)..controls (7,0.75) and (7,-0.75)..(5,-1.5);

\draw (12.3,0) node {$\displaystyle{\ne 0}\quad\Rightarrow\quad x_r\ne x_s , \textrm{ for all } r\ne s$};
\end{tikzpicture}
\end{center}

Moreover, if $\max_i\{x_i\}=a-1$, then the diagram from above can be nonzero if and only if there exists a permutation $\sigma$ of $\{0,1,\ldots,a-1\}$, such that $x_{a-i}=\sigma_i$, $i=0,\ldots,a-1$, in which case  

\begin{center}
\begin{tikzpicture}[scale=0.55]
\draw (2.1,0) node {$\scriptscriptstyle{x_1}$};
\draw (4,0) node {$\scriptscriptstyle{x_2}$};
\draw (5.8,0) node {$\scriptscriptstyle{x_{a-1}}$};
\draw (7.6,0) node {$\scriptstyle{x_{a}}$};
\draw (5,0) node {$\cdots$};
\filldraw[black] 
                  (7.23,0) circle (1.8pt)
                 (6.5,0) circle (1.8pt)
                 
(3.522,0) circle (1.8pt)
                 (2.725,0) circle (1.8pt);


\draw (5.2,2) node {$\scriptstyle{a}$};
\draw (5.2,-2) node {$\scriptstyle{a}$};

\draw [very thick] (5,-2)--(5,-1.5);
\draw [very thick] (5,2)--(5,1.5);
\draw (5,1.5)..controls (2,0.75) and (2,-0.75)..(5,-1.5);
\draw (5,1.5)..controls (8,0.75) and (8,-0.75)..(5,-1.5);
\draw (5,1.5)..controls (3,0.75) and (3,-0.75)..(5,-1.5);
\draw (5,1.5)..controls (7,0.75) and (7,-0.75)..(5,-1.5);

\draw (10,0) node {$\displaystyle{=\sgn \sigma}$};

\draw [very thick] (12,-2)--(12,2);
\draw (12.2,-2.1) node {$\scriptstyle{a}$}; 
\end{tikzpicture}
\end{center}
\end{lemma}

The above lemma implies the following

\begin{lemma}\label{vazlem}\cite[Proposition 2.4.1]{thick}
Let $\alpha\in P(a,x)$ and $\beta\in P(b,y)$ be partitions. Then we have that 
\begin{center}
\begin{tikzpicture}[scale=0.5]
\draw [very thick] (11,-2)--(11,-1.2);
\draw [very thick] (11,2)--(11,1.2);
\draw (11,1.2)..controls (10,0)..(11,-1.2);
\draw (11,1.2)..controls (12,0)..(11,-1.2);
\draw (9.7,0) node {$\scriptstyle{\pi_{\alpha}}$};
\draw (12.3,0) node {$\scriptstyle{\pi_{\beta}}$};
\draw (13.5,0) node {$= s$};
\draw [very thick] (15,-2)--(15,2);

\draw (15.4,-2.2) node {$\scriptstyle{a+b}$};
\draw (15.55,0) node {$\scriptstyle{\pi_{\gamma}}$};
\draw (10.6,-2.3) node {$\scriptstyle{a+b}$};
\draw (10.6,2.3) node {$\scriptstyle{a+b}$};

\draw (10.5,-1) node {$\scriptstyle{a}$};
\draw (11.5,-1) node {$\scriptstyle{b}$};

\filldraw[black] (15,0) circle (2pt);
\filldraw[black] (10.2,0) circle (2pt);
\filldraw[black] (11.8,0) circle (2pt);

\end{tikzpicture}
\end{center}
\noindent for some partition $\gamma\in P(a+b,\max\{x-b,y-a\})$ and $s\in\{-1,0,1\}$. If $s\ne 0$, then $|\gamma|=|\alpha|+|\beta|-ab$.

Moreover, if $\alpha\in P(a,b)$ and $\beta\in P(b,a)$, then we have:
\begin{center}
\begin{tikzpicture}[scale=0.5]
\draw [very thick] (9.5,-2)--(9.5,-1.2);
\draw [very thick] (9.5,2)--(9.5,1.2);
\draw (9.5,1.2)..controls (8.5,0)..(9.5,-1.2);
\draw (9.5,1.2)..controls (10.5,0)..(9.5,-1.2);
\draw (8.2,0) node {$\scriptstyle{\pi_{\alpha}}$};
\draw (10.8,0) node {$\scriptstyle{\pi_{\beta}}$};
\draw (13.5,0) node {$= \delta_{\beta,\hat{\alpha}} (-1)^{|\beta|}$};
\draw [very thick] (16.2,-2)--(16.2,2);

\draw (16.4,-2.2) node {$\scriptstyle{a+b}$};
\draw (10,-2.3) node {$\scriptstyle{a+b}$};
\draw (10,2.3) node {$\scriptstyle{a+b}$};

\draw (9,-1) node {$\scriptstyle{a}$};
\draw (10,-1) node {$\scriptstyle{b}$};

\filldraw[black] (8.7,0) circle (2pt);
\filldraw[black] (10.3,0) circle (2pt);

\end{tikzpicture}
\end{center}
In particular, the left hand side can be nonzero only when $|\alpha|+|\beta|=ab$.
\end{lemma}

\begin{lemma}\cite{canbas}\label{pomocna}
Let $\gamma\in P(a)$ and $\psi\in P(a,b)$ be partitions. Then
\begin{center}
\begin{tikzpicture}[scale=0.5]
\draw [very thick] (11,-2)--(11,-1.2);
\draw [very thick] (11,2)--(11,1.2);
\draw (11,1.2)..controls (10,0)..(11,-1.2);
\draw (11,1.2)..controls (12,0)..(11,-1.2);
\draw (9.7,-0.4) node {$\scriptstyle{\pi_{\gamma}}$};
\draw (9.7,0.4) node {$\scriptstyle{\pi_{\psi}}$};
\draw (13.7,0) node {$=$};
\draw [very thick] (15,-2)--(15,2);

\draw (15.4,-2.2) node {$\scriptstyle{a+b}$};
\draw (17.0,-0.2) node {$\scriptstyle{\pi_{_{^{\gamma}/_{(K_{a,b}-\psi)}}}}$};
\draw (10.6,-2.3) node {$\scriptstyle{a+b}$};
\draw (10.6,2.3) node {$\scriptstyle{a+b}$};

\draw (10.5,-1) node {$\scriptstyle{a}$};
\draw (11.5,-1) node {$\scriptstyle{b}$};

\filldraw[black] (15,0) circle (2pt);
\filldraw[black] (10.33,-0.3) circle (2pt);
\filldraw[black] (10.33,0.3) circle (2pt);


\end{tikzpicture}
\end{center}
\end{lemma}

So far we have been examining the diagrams of a single color. We use the following convention when drawing the diagrams involving two adjacent colors.



{\bf {Notation convention:}} For two colors (indices) that satisfy $i\cdot j=-1$, we shall draw strands colored $i$ as straight lines, and strands colored $j$ as curly lines:
\begin{center}
\begin{tikzpicture}[scale=0.55]

\draw (0,-0.2) node {$\displaystyle{\Id\nolimits_{\mathcal{E}_i^{(a)}}}:$};

\draw [thick] (2,-1)--(2,1);

\draw (1.8,-1.1) node {$\scriptstyle{a}$};

\draw (8,-0.2) node {$\displaystyle{\Id\nolimits_{\mathcal{E}_j^{(b)}}}:$};

\draw [snake=snake,segment amplitude=.2mm,segment length=2mm,thick]
(10,-1)--(10,1);
 
 \draw (9.8,-1.1) node {$\scriptstyle{b}$};
\end{tikzpicture}
\end{center}

Thus, from now on, each line carries one label, and that label represents the thickness of a line.\\
\\

The first ''thick" property that we shall frequently use is about sliding the thick dots past crossings which involve strands of different colors. It follows straightforward from the analogous property for thin strands (\ref{thinslide}), the definition of thick dot (\ref{decdef}) and  associativity of splitters. 

\vskip 0.5cm

\begin{proposition}[Dot Slide]
The thick dots can be freely moved through the thick crossing of the two thick strands with different colors, i.e.:

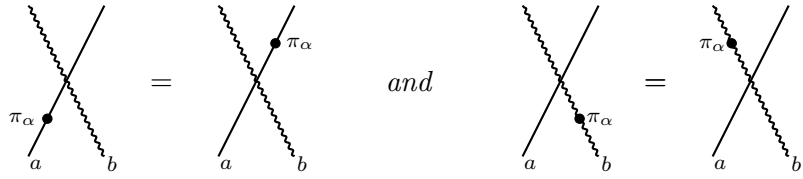
\begin{center}
\begin{tikzpicture} [scale=0.5]

\draw [thick] (-1,-2)--(1,2);
\draw [thick,snake=snake,segment amplitude=.2mm,segment length=1mm] (1,-2)--(-1,2);

\draw (-0.8,-2.2) node {$\scriptstyle{a}$};
\draw (1.2,-2.2) node {$\scriptstyle{b}$};

\draw (2.5,-0.1) node {$=$};

\draw [thick] (4,-2)--(6,2);
\draw [thick,snake=snake,segment amplitude=.2mm,segment length=1mm] (6,-2)--(4,2);

\draw (6.15,1) node {$\scriptstyle{\pi_{\alpha}}$};
\draw (-1.15,-1) node {$\scriptstyle{\pi_{\alpha}}$};

\filldraw[black] (-0.5,-1) circle (3.5pt) 
(5.5,1) circle (3.5pt);

\draw (4.2,-2.2) node {$\scriptstyle{a}$};
\draw (6.2,-2.2) node {$\scriptstyle{b}$};

\draw (9,0) node {{\textrm{and}}};

\draw [thick] (12,-2)--(14,2);
\draw [thick,snake=snake,segment amplitude=.2mm,segment length=1mm] (14,-2)--(12,2);

\draw (12.2,-2.2) node {$\scriptstyle{a}$};
\draw (14.2,-2.2) node {$\scriptstyle{b}$};

\draw (15.5,-0.1) node {$=$};

\draw [thick] (17,-2)--(19,2);
\draw [thick,snake=snake,segment amplitude=.2mm,segment length=1mm] (19,-2)--(17,2);

\draw (17,0.95) node {$\scriptstyle{\pi_{\alpha}}$};
\draw (14.05,-1) node {$\scriptstyle{\pi_{\alpha}}$};

\filldraw[black] (13.5,-1) circle (3.5pt) 
(17.5,1) circle (3.5pt);

\draw (15.5,-0.1) node {$=$};

\draw (17.2,-2.2) node {$\scriptstyle{a}$};
\draw (19.2,-2.2) node {$\scriptstyle{b}$};

\end{tikzpicture}
\end{center}

\end{proposition}

\vskip 0.5cm

The following two propositions  are extensions of the thin R2 and R3 relations (\ref{r2tan}) and (\ref{r3tan}). 

\begin{proposition}[Thick R2 Move]\label{pr2}
We have

\begin{center}
\begin{tikzpicture} [scale=0.5]
\draw [thick] (-1,-3)--(1,0);
\draw [thick] (1,0)--(-1,3);
\draw [thick,snake=snake,segment amplitude=.2mm,segment length=2mm](1,-3)--(-1,0);
\draw [thick,snake=snake,segment amplitude=.2mm,segment length=2mm](-1,0)--(1,3);
\draw (-1.3,-3) node {$\scriptstyle{a}$};
\draw (1.3,-3) node {$\scriptstyle{b}$};
\draw [thick] (3.7,0) node {$= \displaystyle{\sum_{\alpha\in P(a,b)}}$};
\draw [thick](6.5,-3)--(6.5,3);
\draw [thick,snake=snake,segment amplitude=.2mm,segment length=2mm](8,-3)--(8,3);
\draw (5.9,0) node {$\scriptstyle{\pi_{{}_{\alpha}}}$};
\draw (8.55,0) node {$\scriptstyle{\pi_{{}_{\widehat{\alpha}}}}$};
\draw (6.3,-3.2) node {$\scriptstyle{a}$};
\draw (8.1,-3.2) node {$\scriptstyle{b}$};
\filldraw[black] (6.5,0) circle (2pt) 
(8,0) circle (2pt);
\end{tikzpicture}
\end{center}
\end{proposition}

\begin{proposition}[Thick R3 Move]\label{pr3}
We have

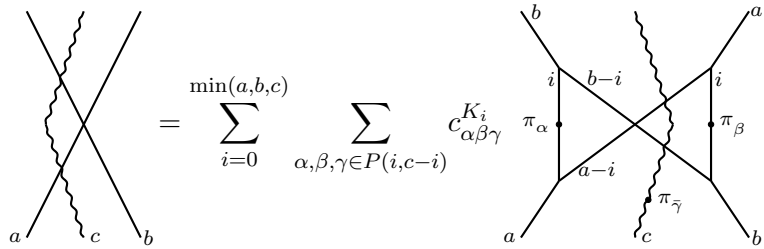
\begin{center}
\begin{tikzpicture} [scale=0.5]
\draw [thick] (-1,-3)--(2,3);
\draw [thick](2,-3)--(-1,3);
\draw [thick, snake=snake,segment amplitude=.2mm,segment length=2mm](0.5,-3)--(-0.5,0);
\draw [thick, snake=snake,segment amplitude=.2mm,segment length=2mm](-0.5,0)--(0.5,3);

\draw (-1.3,-3) node {$\scriptstyle{a}$};
\draw (0.8,-3) node {$\scriptstyle{c}$};
\draw (2.2,-3) node {$\scriptstyle{b}$};

\draw (7,0) node {$= \displaystyle{\sum_{i=0}^{\min(a,b,c)}\sum_{\alpha,\beta,\gamma\in P(i,c-i)}}c^{K_i}_{\alpha\beta\gamma}$};
\draw [thick](12,-3)--(13,-1.5);
\draw [thick](12,3)--(13,1.5);
\draw [thick](18,-3)--(17,-1.5);
\draw [thick](18,3)--(17,1.5);

\draw [thick](13,1.5)--(17,-1.5);
\draw [thick](13,-1.5)--(17,1.5);
\draw [thick](13,1.5)--(13,-1.5);
\draw [thick](17,1.5)--(17,-1.5);

\draw [thick,snake=snake,segment amplitude=.2mm,segment length=2mm](15,-3)--(16,0);
\draw [thick,snake=snake,segment amplitude=.2mm,segment length=2mm](16,0)--(15,3);

\draw (12.4,3) node {$\scriptstyle{b}$};
\draw (18.2,3) node {$\scriptstyle{a}$};
\draw (11.7,-3) node {$\scriptstyle{a}$};
\draw (18.2,-3) node {$\scriptstyle{b}$};
\draw (15.3,-3) node {$\scriptstyle{c}$};

\draw (12.8,1.2) node {$\scriptstyle{i}$};
\draw (17.2,1.2) node {$\scriptstyle{i}$};

\draw (14.2,1.2) node {$\scriptstyle{b-i}$};
\draw (14,-1.2) node {$\scriptstyle{a-i}$};

\draw (12.42,0) node {$\scriptstyle{\pi_{{}_{\alpha}}}$};
\draw (17.6,-0.05) node {$\scriptstyle{\pi_{{}_{\beta}}}$};
\draw (15.96,-2.1) node {$\scriptstyle{\pi_{{}_{\bar{\gamma}}}}$};

\filldraw[black] (13,0) circle (2pt) 
(17,0) circle (2pt)
(15.35,-2) circle (2pt);
\end{tikzpicture}
\end{center}
where $K_i=(\underbrace{c-i,c-i,\ldots,c-i}_i)$, for $i>0$, and $K_0=0$.
\end{proposition}

\section{Proof of Proposition \ref{pr2} - Thick R2 move}\label{spr2}
First we prove the case $a=1$. We rewrite the left hand side by ''exploding" the thick edge into $b$ thin strands by using 
(\ref{decdef})  for $\alpha=0$, then we use the pitchfork lemma and dot slide, and then apply thin R2 move (\ref{r2tan}) $b$ times:
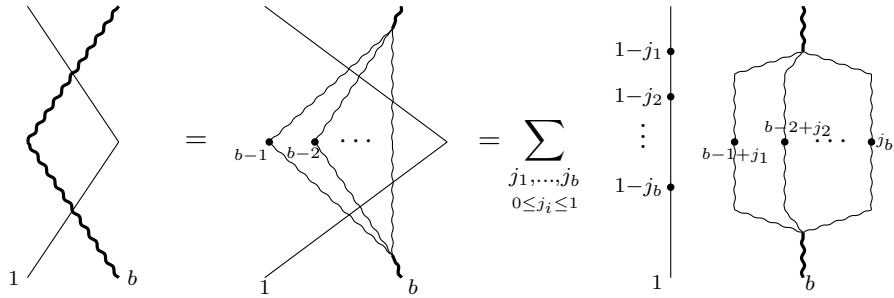
\begin{center}
\begin{tikzpicture} [scale=0.6]
\draw  (-1,-3)--(1,0);
\draw  (1,0)--(-1,3);
\draw [very thick,snake=snake,segment amplitude=.2mm,segment length=2mm](1,-3)--(-1,0);
\draw [very thick,snake=snake,segment amplitude=.2mm,segment length=2mm](-1,0)--(1,3);
\draw (-1.3,-3) node {$\scriptstyle{1}$};
\draw (1.3,-3) node {$\scriptstyle{b}$};
\draw [thick] (2.7,0) node {$=$};

\draw (4.2,-3)--(8.2,0);
\draw (8.2,0)--(4.2,3);


\draw [very thick,snake=snake,segment amplitude=.2mm,segment length=2mm](7.2,-3)--(7,-2.5);
\draw [very thick,snake=snake,segment amplitude=.2mm,segment length=2mm](7.2,3)--(7,2.5);
\draw (4.2,-3.2) node {$\scriptstyle{1}$};
\draw (7.45,-3.2) node {$\scriptstyle{b}$};

\draw [snake=snake,segment amplitude=.1mm,segment length=2mm] (7,-2.5)--(4.3,0);
\draw [snake=snake,segment amplitude=.1mm,segment length=2mm]  (4.3,0)--(7,2.5);
\draw [snake=snake,segment amplitude=.1mm,segment length=2mm] (7,-2.5)--(5.3,0);
\draw [snake=snake,segment amplitude=.1mm,segment length=2mm]  (5.3,0)--(7,2.5);
\draw [snake=snake,segment amplitude=.1mm,segment length=2mm] (7,-2.5)--(7,0);
\draw [snake=snake,segment amplitude=.1mm,segment length=2mm]  (7,0)--(7,2.5);

\draw (3.9,-0.3) node {$\scriptscriptstyle{b-1}$};
\draw (5.05,-0.25) node {$\scriptscriptstyle{b-2}$};
\draw (6.2,0) node {${\cdots}$};


\filldraw[black] (4.3,0) circle (2pt) 
(5.3,0) circle (2pt);


\draw (10,-0.25) node {$=\displaystyle{\sum\limits_{j_1,\ldots,j_b}}$};
\draw (10.35,-1.4) node {$\scriptscriptstyle{0\le j_i\le 1}$};

\draw (13.1,-3)--(13.1,3);
\draw (12.8,-3.1) node {$\scriptstyle{1}$};

\filldraw[black] (13.1,-1) circle (2pt)
(13.1,1) circle (2pt) 
(13.1,2) circle (2pt);

\draw (12.4,-1) node {$\scriptstyle{1-j_b}$};
\draw (12.4,1) node {$\scriptstyle{1-j_2}$};
\draw (12.4,2) node {$\scriptstyle{1-j_1}$};
\draw (12.5,0.3) node {${\vdots}$};

\draw [very thick,snake=snake,segment amplitude=.2mm,segment length=2mm](16,-3)--(16,-2);
\draw [very thick,snake=snake,segment amplitude=.2mm,segment length=2mm](16,3)--(16,2);

\draw (16.15,-3.1) node {$\scriptstyle{b}$};

\draw [snake=snake,segment amplitude=.1mm,segment length=2mm] (16,-2)--(14.5,-1.5);
\draw [snake=snake,segment amplitude=.1mm,segment length=2mm]  (14.5,1.5)--(16,2);
\draw [snake=snake,segment amplitude=.1mm,segment length=2mm] (14.5,-1.5)--(14.5,1.5);

\draw [snake=snake,segment amplitude=.1mm,segment length=2mm] (16,-2)--(15.6,-1.5);
\draw [snake=snake,segment amplitude=.1mm,segment length=2mm]  (15.6,1.5)--(16,2);
\draw [snake=snake,segment amplitude=.1mm,segment length=2mm] (15.6,-1.5)--(15.6,1.5);

\draw [snake=snake,segment amplitude=.1mm,segment length=2mm] (16,-2)--(17.5,-1.5);
\draw [snake=snake,segment amplitude=.1mm,segment length=2mm]  (17.5,1.5)--(16,2);
\draw [snake=snake,segment amplitude=.1mm,segment length=2mm] (17.5,-1.5)--(17.5,1.5);

\filldraw[black] (14.5,0) circle (2pt)
(15.6,0) circle (2pt) 
(17.5,0) circle (2pt);

\draw (14.53,-0.3) node {$\scriptscriptstyle{b-1+j_1}$};
\draw (15.9,0.3) node {$\scriptscriptstyle{b-2+j_2}$};
\draw (17.8,0) node {$\scriptscriptstyle{j_b}$};

\draw (16.6,0) node {${\cdots}$};

\end{tikzpicture}
\end{center}

Since the "exploded" thick strand is antisymmetric with respect to the exchange of the number of dots on two thin strands, we have that the last diagram is zero whenever we have two strands with the same number of dots. Since the number of dots on $i$-th strand is equal to $b-i+j_i$, $i=1,\ldots,b$, and since all $j_i$'s are either $0$ or $1$, the last diagram can be nonzero if and only if $j_1=j_2=\ldots=j_k=1$ and $j_i=0$, for $i>k$, for some $k=0,\ldots,b$. Therefore, the last sum  becomes simply

\begin{center}
\begin{tikzpicture} [scale=0.35]

\draw [thick] (2.8,0) node {$ \displaystyle{\sum_{l+k=b}}$};
\draw (6.5,-3)--(6.5,3);
\draw [very thick,snake=snake,segment amplitude=.2mm,segment length=2mm](8,-3)--(8,3);
\draw (5.9,0) node {$\scriptstyle{l}$};
\draw (8.7,0) node {$\scriptstyle{\varepsilon_{k}}$};
\draw (6.3,-3.2) node {$\scriptstyle{1}$};
\draw (8.2,-3.2) node {$\scriptstyle{b}$};
\filldraw[black] (6.5,0) circle (3pt) 
(8,0) circle (4pt);
\end{tikzpicture}
\end{center}
where $\varepsilon_k$ is the $k$-th elementary symmetric polynomial, as wanted.\\

As for the general case, now we split the strand of thickness $a$ into $a$ thin lines, and again after using pitchfork lemma and dot slides, we apply the thin R2 move $a$ times on the diagram involving strands of thicknesses $1$ and $b$:

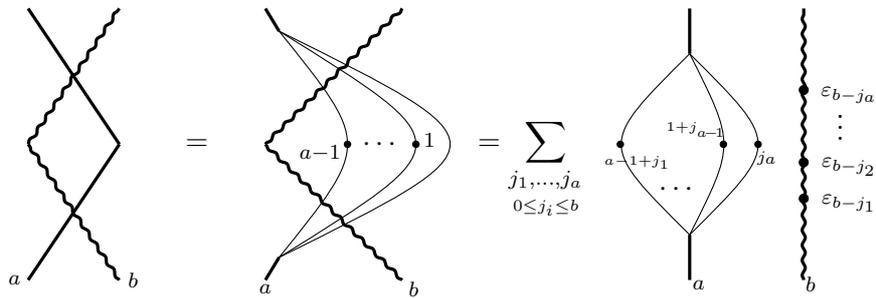
\begin{center}
\begin{tikzpicture} [scale=0.6]
\draw [very thick] (-1,-3)--(1,0);
\draw [very thick] (1,0)--(-1,3);
\draw [very thick,snake=snake,segment amplitude=.2mm,segment length=2mm](1,-3)--(-1,0);
\draw [very thick,snake=snake,segment amplitude=.2mm,segment length=2mm](-1,0)--(1,3);
\draw (-1.3,-3) node {$\scriptstyle{a}$};
\draw (1.3,-3) node {$\scriptstyle{b}$};
\draw [thick] (2.7,0) node {$=$};

\draw [very thick] (4.2,-3)--(4.5,-2.5);
\draw [very thick] (4.5,2.5)--(4.2,3);

\draw (4.5,-2.5)..controls (9.5,0).. (4.5,2.5);
\draw (4.5,-2.5)..controls (8.5,0).. (4.5,2.5);
\draw (4.5,-2.5)..controls (6.5,0).. (4.5,2.5);

\draw [very thick,snake=snake,segment amplitude=.2mm,segment length=2mm](7.2,-3)--(4.2,0);
\draw [very thick,snake=snake,segment amplitude=.2mm,segment length=2mm](7.2,3)--(4.2,0);

\draw (4.2,-3.2) node {$\scriptstyle{a}$};
\draw (7.45,-3.2) node {$\scriptstyle{b}$};

\draw (5.4,-0.2) node {$\scriptstyle{a-1}$};
\draw (7.8,0.1) node {$\scriptstyle{1}$};
\draw (6.7,0) node {${\cdots}$};

\filldraw[black] (6,0) circle (2pt) 
                       (7.5,0) circle (2pt);

\draw (10,-0.25) node {$=\displaystyle{\sum\limits_{j_1,\ldots,j_a}}$};
\draw (10.35,-1.4) node {$\scriptscriptstyle{0\le j_i\le b}$};

\draw [very thick] (13.5,-3)--(13.5,-2);
\draw [very thick] (13.5,3)--(13.5,2);
\draw (13.7,-3.1) node {$\scriptstyle{a}$};

\draw (13.5,-2)..controls (15.5,0).. (13.5,2);
\draw (13.5,-2)..controls (14.5,0).. (13.5,2);
\draw (13.5,-2)..controls (11.5,0).. (13.5,2);

\filldraw[black] (11.98,0) circle (2pt)
                      (14.24,0) circle (2pt) 
                       (15,0) circle (2pt);

\draw (12.35,-0.4) node {\tiny{$\scriptscriptstyle{a-1+j_1}$}};
\draw (13.65,0.3) node {\tiny{$\scriptscriptstyle{1+j_{_{a\!-\!1}}}$}};
\draw (15.2,-0.3) node {$\scriptscriptstyle{j_a}$};
\draw (13.2,-1) node {${\cdots}$};

\filldraw[black] 
                  (16,1.2) circle (2.8pt)
                 (16,-0.4) circle (2.8pt)
                  (16,-1.2) circle (2.8pt);
                 
\draw (17,1.1) node {$\scriptstyle{\varepsilon_{b-j_a}}$};              
\draw (17,-0.5) node {$\scriptstyle{\varepsilon_{b-j_2}}$};
\draw (17,-1.3) node {$\scriptstyle{\varepsilon_{b-j_1}}$};
\draw (16.8,0.55) node {${\vdots}$};

\draw [very thick,snake=snake,segment amplitude=.2mm,segment length=2mm](16,-3)--(16,3);

\draw (16.15,-3.1) node {$\scriptstyle{b}$};

\end{tikzpicture}
\end{center}
 
Since $\varepsilon_k$ as a polynomial in $b$ variables is zero whenever $k<0$ or $k>b$, in the last expression we can take summations over all integers $j_i$. Let $x_i=a-i+j_i$, $i=1,\ldots,a$. Then the last expression becomes 

\begin{center}
\begin{tikzpicture}[scale=0.55]
\draw (0,-0.2) node {$\displaystyle{\sum_{x_1,\ldots,x_a}}$};

\draw (2.1,0) node {$\scriptscriptstyle{x_1}$};
\draw (4,0) node {$\scriptscriptstyle{x_2}$};
\draw (5.85,0) node {$\scriptscriptstyle{x_{a-1}}$};
\draw (7.6,0) node {$\scriptstyle{x_{a}}$};
\draw (4.8,0) node {$\cdots$};
\filldraw[black] 
                  (7.23,0) circle (1.8pt)
                 (6.5,0) circle (1.8pt)
                 
(3.522,0) circle (1.8pt)
                 (2.725,0) circle (1.8pt);


\draw (5.2,-2) node {$\scriptstyle{a}$};

\draw [very thick] (5,-2)--(5,-1.5);
\draw [very thick] (5,2)--(5,1.5);
\draw (5,1.5)..controls (2,0.75) and (2,-0.75)..(5,-1.5);
\draw (5,1.5)..controls (8,0.75) and (8,-0.75)..(5,-1.5);
\draw (5,1.5)..controls (3,0.75) and (3,-0.75)..(5,-1.5);
\draw (5,1.5)..controls (7,0.75) and (7,-0.75)..(5,-1.5);


\draw [very thick,snake=snake,segment amplitude=.2mm,segment length=2mm](9.5,-2)--(9.5,2);
\draw (9.7,-2.1) node {$\scriptstyle{b}$}; 

\filldraw[black] 
                  (9.5,1.2) circle (2.8pt)
                 (9.5,-0.4) circle (2.8pt)
                  (9.5,-1.2) circle (2.8pt);
                 
\draw (11,1.1) node {$\scriptstyle{\varepsilon_{b+a-a-x_a}}$};              
\draw (11,-0.5) node {$\scriptstyle{\varepsilon_{b+a-2-x_2}}$};
\draw (11,-1.3) node {$\scriptstyle{\varepsilon_{b+a-1-x_1}}$};
\draw (10.1,0.55) node {${\vdots}$};

\end{tikzpicture}
\end{center}

Due to antisymmetry of thin lines in the "exploded" thick strand, this can be nonzero only when all $x_i$'s are pairwise distinct. Moreover, if we denote the decreasing ordering of such $x_i$'s by $y_1>y_2>\ldots>y_a$, then for some permutation $\sigma\in S_a$, we have $x_i=y_{\sigma_i}$, $i=1,\ldots,a$, and the whole expression becomes

\begin{center}
\begin{tikzpicture}[scale=0.55]
\draw (-1.2,-0.2) node 
{$\displaystyle{\sum_{y_1>\cdots>y_a} \sum_{\sigma\in S_a} \sgn \sigma}$};

\draw (2.35,0) node {$\scriptscriptstyle{y_1}$};
\draw (3.85,0) node {$\scriptscriptstyle{y_2}$};
\draw (5.8,0) node {$\scriptscriptstyle{y_{a-1}}$};
\draw (7.6,0) node {$\scriptstyle{y_{a}}$};
\draw (4.75,0) node {$\cdots$};
\filldraw[black] 
                  (7.23,0) circle (1.8pt)
                 (6.5,0) circle (1.8pt)
                 (3.522,0) circle (1.8pt)
                 (2.725,0) circle (1.8pt);


\draw (5.2,-2) node {$\scriptstyle{a}$};

\draw [very thick] (5,-2)--(5,-1.5);
\draw [very thick] (5,2)--(5,1.5);
\draw (5,1.5)..controls (2,0.75) and (2,-0.75)..(5,-1.5);
\draw (5,1.5)..controls (8,0.75) and (8,-0.75)..(5,-1.5);
\draw (5,1.5)..controls (3,0.75) and (3,-0.75)..(5,-1.5);
\draw (5,1.5)..controls (7,0.75) and (7,-0.75)..(5,-1.5);


\draw [very thick,snake=snake,segment amplitude=.2mm,segment length=2mm](9,-2)--(9,2);
\draw (9.2,-2.1) node {$\scriptstyle{b}$}; 

\filldraw[black] 
                  (9,1.2) circle (2.8pt)
                 (9,-0.4) circle (2.8pt)
                  (9,-1.2) circle (2.8pt);
                 
\draw (10.7,1.1) node {$\scriptstyle{\varepsilon_{b+a-a-y_{\sigma_a}}}$};              
\draw (10.7,-0.5) node {$\scriptstyle{\varepsilon_{b+a-2-y_{\sigma_2}}}$};
\draw (10.7,-1.3) node {$\scriptstyle{\varepsilon_{b+a-1-y_{\sigma_1}}}$};
\draw (9.8,0.55) node {${\vdots}$};

\end{tikzpicture}
\end{center}

\noindent By setting $\alpha_i=y_i-(a-i)$, $i=1,\ldots,a$, and $\alpha=(\alpha_1,\ldots,\alpha_a)$, by (\ref{decdef}) we have

\begin{center}
\begin{tikzpicture}[scale=0.55]

\draw (14.2,-0.2) node 
{$\displaystyle{\sum_{\alpha\in P(a)} \sum_{\sigma\in S_a} \sgn \sigma}$};

\draw [very thick] (18,-2)--(18,2);
\draw [very thick,snake=snake,segment amplitude=.2mm,segment length=2mm](19.5,-2)--(19.5,2);

\filldraw[black] 
                  (18,0) circle (2.8pt);

\draw (17.5,0) node {$\scriptstyle{\pi_{\alpha}}$};
\draw (18.2,-2.1) node {$\scriptstyle{a}$};

\draw (19.7,-2.1) node {$\scriptstyle{b}$}; 

\filldraw[black] 
                  (19.5,1.2) circle (2.8pt)
                 (19.5,-0.4) circle (2.8pt)
                  (19.5,-1.2) circle (2.8pt);
                 
\draw (21.4,1.1) node {$\scriptstyle{\varepsilon_{b-a-\alpha_{_{\sigma_a}}+\sigma_a}}$};              
\draw (21.4,-0.5) node {$\scriptstyle{\varepsilon_{b-2-\alpha_{_{\sigma_2}}+\sigma_2}}$};
\draw (21.4,-1.3) node {$\scriptstyle{\varepsilon_{b-1-\alpha_{_{\sigma_1}}+\sigma_1}}$};
\draw (20.6,0.55) node {${\vdots}$};

\end{tikzpicture}
\end{center}

\noindent Now, from the definition of the determinant of a matrix we have
\begin{eqnarray*}
&\sum\limits_{\sigma\in S_a} \sgn \sigma \prod_{i=1}^{a} \varepsilon_{b-i-\alpha_{_{\sigma_i}}+\sigma_i}= \det [\varepsilon_{b-i-\alpha_j+j}]_{i,j=1}^{a}=\\
&\quad\quad\quad\quad=\det [\varepsilon_{b-\alpha_{a+1-i}+j-i}]_{i,j=1}^{a}=\pi_{\overline{(b-\alpha_a,\ldots,b-\alpha_2,b-\alpha_1)}}.
\end{eqnarray*}
In the second equality we used the fact that the determinants of a matrix and its transpose are equal, and that the determinant doesn't change if we invert the orders of the rows and of the columns of a matrix. Finally, the last equality is the Giambelli formula (\ref{djam}).  Therefore, in the last diagram, a summand can be nonzero only for $\alpha\in P(a,b)$, and thus we get the desired equality. \kraj

\section{Proof of Proposition \ref{pr3} - Thick R3 move}\label{spr3}

We will prove the proposition by induction on the thicknesses of the three strands involved. We shall frequently use the following 
\begin{equation}\label{pomoc11}
\begin{tikzpicture} [scale=0.5]
\draw [very thick] (0,-2)--(0,2);

\draw (1.5,-0.1) node {$=$};

\draw [very thick] (4,-2)--(4,-1.5);
\draw [very thick] (4,2)--(4,1.5);

\draw (4,-1.5) .. controls (3,0) .. (4,1.5);
\draw [thick] (4,-1.5) .. controls (5,0) .. (4,1.5);

\draw (0.7,-2.1) node {$\scriptstyle{a+1}$};

\draw (4.7,-2.1) node {$\scriptstyle{a+1}$}; 
\draw (4.7,2.1) node {$\scriptstyle{a+1}$};
\draw (3.6,-1.25) node {$\scriptscriptstyle{1}$};
\draw (4.5,-1.25) node {$\scriptscriptstyle{a}$};
\draw (3,0.2) node {$\scriptscriptstyle{a}$};

\filldraw[black] (3.27,0) circle (2.5pt);

\end{tikzpicture}
\end{equation}
which follows directly from (\ref{decdef}) and Associativity of splitters.

First we prove the case when $a=c=1$, i.e. we show that for any $b$:

\begin{equation}\label{11b}
\begin{tikzpicture}[scale=0.5]

\draw (0,-2)--(4,2);
\draw [thick] (4,-2)--(0,2);
\draw [snake=snake,segment amplitude=.2mm,segment length=2mm] (2,-2)--(1,0);
\draw [snake=snake,segment amplitude=.2mm,segment length=2mm] (2,2)--(1,0);

\draw (0.2,-2.1) node {$\scriptstyle{1}$};
\draw (2.2,-2.1) node {$\scriptstyle{1}$};
\draw (4.2,-2.1) node {$\scriptstyle{b}$};

\draw (5.5,-0.1) node {$=$};

\draw (7,-2)--(11,2);
\draw [thick] (11,-2)--(7,2);
\draw [snake=snake,segment amplitude=.2mm,segment length=2mm] (9,-2)--(10,0);
\draw [snake=snake,segment amplitude=.2mm,segment length=2mm] (9,2)--(10,0);

\draw (7.2,-2.1) node {$\scriptstyle{1}$};
\draw (9.2,-2.1) node {$\scriptstyle{1}$};
\draw (11.2,-2.1) node {$\scriptstyle{b}$};

\draw (12.5,-0.1) node {$+$};

\draw (14.5,-2)--(14.5,1);
\draw (17.5,2)--(17.5,-1);

\draw [thick] (14.5,2)--(14.5,1);
\draw [thick] (17.5,-2)--(17.5,-1);

\draw [thick] (17.5,-1)--(14.5,1);
\draw [snake=snake,segment amplitude=.2mm,segment length=2mm] (16,2)--(16,-2);

\draw (14.7,-2.1) node {$\scriptstyle{1}$};
\draw (16.2,-2.1) node {$\scriptstyle{1}$};
\draw (17.7,-2.1) node {$\scriptstyle{b}$};
\draw (16.6,-0.8) node {$\scriptscriptstyle{b-1}$};
\draw (14.7,2.1) node {$\scriptstyle{b}$};
\draw (17.7,2.1) node {$\scriptstyle{1}$};

\end{tikzpicture}
\end{equation}

\noindent For $b=1$ this is simply the thin R3 relation (\ref{r3tan}). Now, suppose that (\ref{11b}) is valid for some $b\ge 1$, and we shall prove that it holds for $b+1$, as well.
We start by using (\ref{pomoc11}) on the strand of thickness $b+1$:
\begin{center}
\begin{tikzpicture}[scale=0.6]

\draw (0,-2)--(3,2);
\draw [thick] (3,-2)--(0,2);
\draw [snake=snake,segment amplitude=.2mm,segment length=2mm] (1.5,-2)--(0.7,0);
\draw [snake=snake,segment amplitude=.2mm,segment length=2mm] (1.5,2)--(0.7,0);

\draw (0.2,-2.1) node {$\scriptstyle{1}$};
\draw (1.7,-2.1) node {$\scriptstyle{1}$};
\draw (3.4,-2.2) node {$\scriptstyle{b+1}$};

\draw (4,-0.1) node {$=$};

\draw (5,-2)--(8,2);
\draw [thick] (8,-2)--(7.6,-1.5);
\draw [thick] (5.4,1.5)--(5,2);
\draw [thick] (5.4,1.5)..controls (7.2,0.7)..(7.6,-1.5);
\draw  (5.4,1.5)..controls (5.8,0)..(7.6,-1.5);

\draw [snake=snake,segment amplitude=.2mm,segment length=2mm] (6.5,-2)--(5.7,0);
\draw [snake=snake,segment amplitude=.2mm,segment length=2mm] (6.5,2)--(5.7,0);

\draw (5.2,-2.1) node {$\scriptstyle{1}$};
\draw (6.7,-2.1) node {$\scriptstyle{1}$};
\draw (8.4,-2.2) node {$\scriptstyle{b+1}$};
\draw (4.7,2.2) node {$\scriptstyle{b+1}$};
\draw (5.35,1.15) node {$\scriptscriptstyle{1}$};
\draw (5.9,1.55) node {$\scriptscriptstyle{b}$};
\draw (6.82,-1.1) node {$\scriptscriptstyle{b}$};

\filldraw[black] 
                  (6.9,-0.9) circle (1.8pt);

\draw (9,-0.1) node {$=$};

\draw (10,-2)--(13,2);
\draw [thick] (13,-2)--(12.6,-1.5);
\draw [thick] (10.4,1.5)--(10,2);
\draw [thick] (10.4,1.5)..controls (12.2,0.7)..(12.6,-1.5);
\draw  (10.4,1.5)..controls (10.8,0)..(12.6,-1.5);

\draw [snake=snake,segment amplitude=.2mm,segment length=2mm] (11.65,-2)--(11.65,2);

\draw (10.2,-2.1) node {$\scriptstyle{1}$};
\draw (11.8,-2.1) node {$\scriptstyle{1}$};
\draw (13.4,-2.2) node {$\scriptstyle{b+1}$};
\draw (9.7,2.2) node {$\scriptstyle{b+1}$};
\draw (10.35,1.15) node {$\scriptscriptstyle{1}$};
\draw (10.9,1.55) node {$\scriptscriptstyle{b}$};
\draw (11.9,-1.2) node {$\scriptscriptstyle{b}$};

\filldraw[black] 
                  (12,-1) circle (1.8pt);

\draw (14,-0.1) node {$+$};

\draw (15,-2)--(15.5,1.2);
\draw (18.5,2)--(17.5,0.5);
\draw [thick] (15.5,1.2).. controls (16.5,1)..(17.5,0.5);

\draw [thick] (15,2)--(15.5,1.2);
\draw [thick] (18.5,-2)--(18.1,-1.5);

\draw [thick] (17.5,0.5)..controls (18,-0.2)..(18.1,-1.5);
\draw  (17.5,0.5)..controls (17,-0.3)..(18.1,-1.5);

\draw [snake=snake,segment amplitude=.2mm,segment length=2mm] (16.5,2)--(16.5,-2);

\draw (14.9,-2.1) node {$\scriptstyle{1}$};
\draw (16.65,-2.1) node {$\scriptstyle{1}$};
\draw (18.9,-2.2) node {$\scriptstyle{b+1}$};
\draw (17.3,-1) node {$\scriptscriptstyle{b}$};
\draw (15.8,1.4) node {$\scriptscriptstyle{b}$};
\draw (14.6,2.1) node {$\scriptstyle{b+1}$};
\draw (18.35,2.1) node {$\scriptstyle{1}$};

\filldraw[black] 
                  (17.65,-1) circle (1.8pt);

\draw (19.5,-0.1) node {$=$};

\end{tikzpicture}
\end{center}

\begin{center}
\begin{tikzpicture}[scale=0.6]

\draw (9,-0.1) node {$=$};

\draw (10,-2)--(12.5,2);
\draw [thick] (13,-2)--(12.6,-1.5);
\draw [thick] (10.4,1.5)--(10,2);
\draw [thick] (10.4,1.5)..controls (12.2,0.7)..(12.6,-1.5);
\draw  (10.4,1.5)..controls (10.8,0)..(12.6,-1.5);

\draw [snake=snake,segment amplitude=.2mm,segment length=2mm] (11.3,-2)--(12.2,0.8);
\draw [snake=snake,segment amplitude=.2mm,segment length=2mm] (11.3,2)--(12.2,0.8);

\draw (10.2,-2.1) node {$\scriptstyle{1}$};
\draw (11.45,-2.1) node {$\scriptstyle{1}$};
\draw (13.4,-2.2) node {$\scriptstyle{b+1}$};
\draw (9.7,2.2) node {$\scriptstyle{b+1}$};
\draw (10.35,1.15) node {$\scriptscriptstyle{1}$};
\draw (10.9,1.55) node {$\scriptscriptstyle{b}$};
\draw (11.9,-1.2) node {$\scriptscriptstyle{b}$};

\filldraw[black] 
                  (12,-1) circle (1.8pt);

\draw (14,-0.1) node {$+$};

\draw (15,-2)--(16.5,0.5);
\draw (18.5,2)--(18,-0.5);
\draw [thick] (15.5,1.2).. controls (16.2,1)..(16.5,0.5);

\draw [thick] (15,2)--(15.5,1.2);
\draw [thick] (18.5,-2)--(18.1,-1.5);
\draw [thick] (18,-0.5)--(16.5,0.5);

\draw [thick] (18,-0.5)--(18.1,-1.5);
\draw  (18.1,-1.5)..controls (15.5,-1)..(15.5,1.2);

\draw [snake=snake,segment amplitude=.2mm,segment length=2mm] (17,2)--(17,-2);

\draw (14.9,-2.1) node {$\scriptstyle{1}$};
\draw (17.15,-2.1) node {$\scriptstyle{1}$};
\draw (18.9,-2.2) node {$\scriptstyle{b+1}$};
\draw (17.5,-1.2) node {$\scriptscriptstyle{b}$};
\draw (15.8,1.37) node {$\scriptscriptstyle{b}$};
\draw (14.6,2.1) node {$\scriptstyle{b+1}$};
\draw (18.35,2.1) node {$\scriptstyle{1}$};
\draw (15.35,0.6) node {$\scriptscriptstyle{1}$};
\draw (18.25,-0.8) node {$\scriptscriptstyle{b}$};
\draw (17.6,0.1) node {$\scriptscriptstyle{b-1}$};

\filldraw[black] 
                  (17.65,-1.41) circle (1.8pt);

\draw (19.5,-0.1) node {$+$};

\draw (20.5,-2)--(20.5,1);
\draw (23.5,2)--(23.5,-1);

\draw [thick] (20.5,2)--(20.5,1);
\draw [thick] (23.5,-2)--(23.5,-1);

\draw [thick] (23.5,-1)--(20.5,1);
\draw [snake=snake,segment amplitude=.2mm,segment length=2mm] (22,2)--(22,-2);

\draw (20.7,-2.1) node {$\scriptstyle{1}$};
\draw (22.2,-2.1) node {$\scriptstyle{1}$};
\draw (24,-2.1) node {$\scriptstyle{b+1}$};
\draw (22.8,-0.8) node {$\scriptscriptstyle{b}$};
\draw (21,2.1) node {$\scriptstyle{b+1}$};
\draw (23.7,2.1) node {$\scriptstyle{1}$};

\end{tikzpicture}
\end{center}

\noindent In the first equality we have used the pitchfork lemma, while the second one follows from the thin R3 move. The third equality follows by applying the induction hypothesis on the first diagram, and by using (\ref{pomoc11}) on the second diagram. The first diagram is exactly the wanted first term on the right hand side of (\ref{11b}), by pitchfork lemma and (\ref{pomoc11}). Finally, the middle diagram in the last expression after applying the associativity of splitters becomes
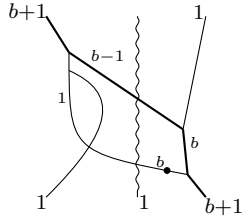
\begin{center}
\begin{tikzpicture}[scale=0.6]

\draw (15,-2).. controls (16.5,-0.5) and (16.6,0.3).. (15.5,0.8);
\draw (18.5,2)--(18,-0.5);

\draw [thick] (15,2)--(15.5,1.2);
\draw [thick] (18.5,-2)--(18.1,-1.5);
\draw [thick] (18,-0.5)--(15.5,1.2);

\draw [thick] (18,-0.5)--(18.1,-1.5);
\draw  (18.1,-1.5)..controls (15.5,-1)..(15.5,1.2);

\draw [snake=snake,segment amplitude=.2mm,segment length=2mm] (17,2)--(17,-2);

\draw (14.9,-2.1) node {$\scriptstyle{1}$};
\draw (17.15,-2.1) node {$\scriptstyle{1}$};
\draw (18.9,-2.2) node {$\scriptstyle{b+1}$};
\draw (17.5,-1.2) node {$\scriptscriptstyle{b}$};
\draw (16.35,1.1) node {$\scriptscriptstyle{b-1}$};
\draw (14.6,2.1) node {$\scriptstyle{b+1}$};
\draw (18.35,2.1) node {$\scriptstyle{1}$};
\draw (15.35,0.2) node {$\scriptscriptstyle{1}$};
\draw (18.25,-0.8) node {$\scriptscriptstyle{b}$};

\filldraw[black] 
                  (17.65,-1.41) circle (1.8pt);

\end{tikzpicture}
\end{center}
\noindent and is equal to zero since it contains a dotless digon (Lemma 2).\\

Now, we pass to the case $c=1$, i.e. we shall prove

\begin{equation}\label{a1b}
\begin{tikzpicture}[scale=0.5]

\draw [thick] (0,-2)--(4,2);
\draw [thick] (4,-2)--(0,2);
\draw [snake=snake,segment amplitude=.2mm,segment length=2mm] (2,-2)--(1,0);
\draw [snake=snake,segment amplitude=.2mm,segment length=2mm] (2,2)--(1,0);

\draw (0.2,-2.1) node {$\scriptstyle{a}$};
\draw (2.2,-2.1) node {$\scriptstyle{1}$};
\draw (4.2,-2.1) node {$\scriptstyle{b}$};

\draw (5.5,-0.1) node {$=$};

\draw [thick] (7,-2)--(11,2);
\draw [thick] (11,-2)--(7,2);
\draw [snake=snake,segment amplitude=.2mm,segment length=2mm] (9,-2)--(10,0);
\draw [snake=snake,segment amplitude=.2mm,segment length=2mm] (9,2)--(10,0);

\draw (7.2,-2.1) node {$\scriptstyle{a}$};
\draw (9.2,-2.1) node {$\scriptstyle{1}$};
\draw (11.2,-2.1) node {$\scriptstyle{b}$};

\draw (12.5,-0.1) node {$+$};

\draw [thick] (14,-2)--(18,2);
\draw [thick] (18,-2)--(14,2);
\draw (14.8,-1.2)--(14.8,1.2);
\draw (17.2,-1.2)--(17.2,1.2);
\draw [snake=snake,segment amplitude=.2mm,segment length=2mm] (16,-2)--(16.7,0);
\draw [snake=snake,segment amplitude=.2mm,segment length=2mm] (16,2)--(16.7,0);

\draw (14.2,-2.1) node {$\scriptstyle{a}$};
\draw (16.2,-2.1) node {$\scriptstyle{1}$};
\draw (18.2,-2.1) node {$\scriptstyle{b}$};
\draw (14.6,-0.1) node {$\scriptscriptstyle{1}$};
\draw (17.4,-0.1) node {$\scriptscriptstyle{1}$};
\draw (14.2,2.1) node {$\scriptstyle{b}$};
\draw (18.2,2.1) node {$\scriptstyle{a}$};
\end{tikzpicture}
\end{equation}
We prove this formula by induction on $a$. We assume that (\ref{a1b}) is valid for some $a\ge 1$ and prove that it also holds for $a+1$. We rewrite the strand of thickness $a+1$ by using (\ref{pomoc11}), and after performing dot slide and pitchfork lemma, we get

\begin{center}
\begin{tikzpicture}[scale=0.6]

\draw [thick] (0,-2)--(3,2);
\draw [thick] (3,-2)--(0,2);
\draw [snake=snake,segment amplitude=.2mm,segment length=2mm] (1.5,-2)--(0.7,0);
\draw [snake=snake,segment amplitude=.2mm,segment length=2mm] (1.5,2)--(0.7,0);

\draw (0.4,-2.1) node {$\scriptstyle{a+1}$};
\draw (1.7,-2.1) node {$\scriptstyle{1}$};
\draw (3.15,-2.2) node {$\scriptstyle{b}$};

\draw (4,-0.1) node {$=$};

\draw [thick] (8,-2)--(5,2);
\draw [thick] (5,-2)--(5.4,-1.5);
\draw [thick] (7.6,1.5)--(8,2);
\draw  (7.6,1.5)..controls (5.8,0.7)..(5.4,-1.5);
\draw [thick] (7.6,1.5)..controls (7.2,0)..(5.4,-1.5);

\draw [snake=snake,segment amplitude=.2mm,segment length=2mm] (7,-2)--(5.5,0.8);
\draw [snake=snake,segment amplitude=.2mm,segment length=2mm] (6.5,2)--(5.5,0.8);

\draw (5.4,-2.2) node {$\scriptstyle{a+1}$};
\draw (6.85,-2.2) node {$\scriptstyle{1}$};
\draw (8.15,-2.2) node {$\scriptstyle{b}$};
\draw (8.2,2.2) node {$\scriptstyle{a+1}$};
\draw (7.1,1.5) node {$\scriptscriptstyle{1}$};
\draw (7.65,1) node {$\scriptscriptstyle{a}$};
\draw (5.3,-1) node {$\scriptscriptstyle{a}$};

\filldraw[black] 
                  (5.5,-0.9) circle (1.8pt);

\draw (9,-0.1) node {$=$};

\draw [thick] (13,-2)--(10,2);
\draw [thick] (10,-2)--(10.4,-1.5);
\draw [thick] (12.6,1.5)--(13,2);
\draw  (12.6,1.5)..controls (10.8,0.7)..(10.4,-1.5);
\draw [thick] (12.6,1.5)..controls (12.2,0)..(10.4,-1.5);

\draw [snake=snake,segment amplitude=.2mm,segment length=2mm] (11.45,-2)--(11.45,2);

\draw (10.2,-2.2) node {$\scriptstyle{a+1}$};
\draw (11.6,-2.2) node {$\scriptstyle{1}$};
\draw (13.15,-2.2) node {$\scriptstyle{b}$};
\draw (13.2,2.2) node {$\scriptstyle{a+1}$};
\draw (12.1,1.5) node {$\scriptscriptstyle{1}$};
\draw (12.65,1) node {$\scriptscriptstyle{a}$};
\draw (10.3,-1) node {$\scriptscriptstyle{a}$};

\filldraw[black] 
                  (10.5,-0.9) circle (1.8pt);

\draw (14,-0.1) node {$+$};

\draw [thick] (15.4,1)--(15,2);
\draw [thick] (18,-2)--(17,0);
\draw [thick] (15.4,1)--(17,0);

\draw [thick] (15,-2)--(15.4,-1.5);
\draw [thick] (17.6,1.5)--(18,2);
\draw  (17.6,1.5)--(17,0);
\draw  (15.4,1)--(15.4,-1.5);

\draw [thick] (17.6,1.5)..controls (18,-1)..(15.4,-1.5);

\draw [snake=snake,segment amplitude=.2mm,segment length=2mm] (16.45,-2)--(16.45,2);

\draw (15.2,-2.2) node {$\scriptstyle{a+1}$};
\draw (16.6,-2.2) node {$\scriptstyle{1}$};
\draw (18.15,-2.2) node {$\scriptstyle{b}$};
\draw (15.15,2.1) node {$\scriptstyle{b}$};
\draw (18.2,2.2) node {$\scriptstyle{a+1}$};
\draw (17.2,1.1) node {$\scriptscriptstyle{1}$};
\draw (17.9,0.9) node {$\scriptscriptstyle{a}$};
\draw (15.2,-1) node {$\scriptscriptstyle{a}$};
\draw (15.2,0.5) node {$\scriptscriptstyle{1}$};
\draw (16,0.94) node {$\scriptscriptstyle{b-1}$};

\filldraw[black] 
                  (15.4,-0.9) circle (1.8pt);
                  
\draw (19.5,-0.1) node {$=$};

\end{tikzpicture}
\end{center}

\begin{center}
\begin{tikzpicture}[scale=0.6]

\draw (4,-0.1) node {$=$};

\draw [thick] (8,-2)--(5,2);
\draw [thick] (5,-2)--(5.4,-1.5);
\draw [thick] (7.6,1.5)--(8,2);
\draw  (7.6,1.5)..controls (5.8,0.7)..(5.4,-1.5);
\draw [thick] (7.6,1.5)..controls (7.2,0)..(5.4,-1.5);

\draw [snake=snake,segment amplitude=.2mm,segment length=2mm] (6.3,2)--(7.5,-0.8);
\draw [snake=snake,segment amplitude=.2mm,segment length=2mm] (7,-2)--(7.5,-0.8);

\draw (5.4,-2.2) node {$\scriptstyle{a+1}$};
\draw (6.85,-2.2) node {$\scriptstyle{1}$};
\draw (8.15,-2.2) node {$\scriptstyle{b}$};
\draw (8.2,2.2) node {$\scriptstyle{a+1}$};
\draw (7.1,1.5) node {$\scriptscriptstyle{1}$};
\draw (7.65,1) node {$\scriptscriptstyle{a}$};
\draw (5.3,-1) node {$\scriptscriptstyle{a}$};

\filldraw[black] 
                  (5.5,-0.9) circle (1.8pt);

\draw (9,-0.1) node {$+$};

\draw [thick] (13,-2)--(10,2);
\draw [thick] (10,-2)--(10.4,-1.5);
\draw [thick] (12.6,1.5)--(13,2);
\draw  (12.6,1.5)..controls (10,1)..(10.4,-1.5);
\draw [thick] (12.6,1.5)..controls (12.2,0)..(10.4,-1.5);

\draw (12.25,-1)--(12.25,0.42);
\draw (11.125,0.5)--(11.125,-0.9);

\draw [snake=snake,segment amplitude=.2mm,segment length=2mm] (12,-2)--(12,2);

\draw (10.2,-2.2) node {$\scriptstyle{a+1}$};
\draw (9.9,2.2) node {$\scriptstyle{b}$};
\draw (12.1,-2.2) node {$\scriptstyle{1}$};
\draw (13.15,-2.2) node {$\scriptstyle{b}$};
\draw (13.2,2.2) node {$\scriptstyle{a+1}$};
\draw (10.1,0.6) node {$\scriptscriptstyle{1}$};
\draw (12.65,1) node {$\scriptscriptstyle{a}$};
\draw (10.1,-1) node {$\scriptscriptstyle{a}$};
\draw (11,-0.2) node {$\scriptscriptstyle{1}$};
\draw (12.45,-0.3) node {$\scriptscriptstyle{1}$};

\filldraw[black] 
                  (10.3,-0.9) circle (1.8pt);

\draw (14,-0.1) node {$+$};

\draw [thick] (15.4,1)--(15,2);
\draw [thick] (18,-2)--(17.6,-1.5);
\draw [thick] (15.4,1)--(17.6,-1.5);

\draw [thick] (15,-2)--(15.4,-1.5);
\draw [thick] (17.6,1.5)--(18,2);
\draw  (17.6,1.5)..controls (16.5,0.4)..(17.6,-1.5);
\draw  (15.4,1)--(15.4,-1.5);

\draw [thick] (17.6,1.5)..controls (17.3,-0.8)..(15.4,-1.5);

\draw [snake=snake,segment amplitude=.2mm,segment length=2mm] (16.45,-2)--(16.45,2);

\draw (15.2,-2.2) node {$\scriptstyle{a+1}$};
\draw (16.6,-2.2) node {$\scriptstyle{1}$};
\draw (18.15,-2.2) node {$\scriptstyle{b}$};
\draw (15.15,2.1) node {$\scriptstyle{b}$};
\draw (18.2,2.2) node {$\scriptstyle{a+1}$};
\draw (17.1,1.25) node {$\scriptscriptstyle{1}$};
\draw (17.75,0.9) node {$\scriptscriptstyle{a}$};
\draw (15.2,-1) node {$\scriptscriptstyle{a}$};
\draw (15.2,0.5) node {$\scriptscriptstyle{1}$};
\draw (16,0.9) node {$\scriptscriptstyle{b-1}$};

\filldraw[black] 
                  (15.4,-0.9) circle (1.8pt);
                  

\end{tikzpicture}
\end{center}

\noindent The second equality follows from (\ref{11b}). The third equality follows by applying the induction hypothesis on the first diagram, and by performing pitchfork lemma on the second diagram.
The third of the three obtained summands is equal to zero since it contains a dotless digon, while the first summand equals the wanted first term from (\ref{a1b}) by using pitchfork lemma and (\ref{pomoc11}). The remaining (second) summand, after applying associativity of splitters and pitchfork lemma becomes:
\begin{center}
\begin{tikzpicture}[scale=0.6]

\draw[thick] (0,-2)--(0.5,-1.5);
\draw[thick] (4,2)--(3.5,1.7);

\draw[thick] (2.5,1.3)--(3.5,1.7);
\draw (3.5,1.7)--(3.5,-1.5);
\draw (1.8,-1)--(1.8,0.2);
\draw[thick] (0,2)--(3.5,-1.5);
\draw[thick] (3.5,-1.5)--(4,-2);

\draw (0.5,-1.5)..controls (1,0.5).. (2.5,1.3);
\draw [thick] (0.5,-1.5)..controls (1.4,-1.4).. (1.8,-1);
\draw [thick] (1.8,-1)..controls (2.5,0.5).. (2.5,1.3);

\draw [snake=snake,segment amplitude=.2mm,segment length=2mm] (3,-2)--(3,2);

\filldraw[black] (0.68,-0.8) circle (1.8pt);

\draw (0.2,-2.2) node {$\scriptstyle{a+1}$};
\draw (-0.1,2.2) node {$\scriptstyle{b}$};
\draw (2.9,-2.2) node {$\scriptstyle{1}$};
\draw (4.15,-2.2) node {$\scriptstyle{b}$};
\draw (4.2,2.2) node {$\scriptstyle{a+1}$};
\draw (0.8,0.35) node {$\scriptscriptstyle{1}$};
\draw (2.72,1.55) node {$\scriptscriptstyle{a}$};
\draw (0.35,-0.8) node {$\scriptscriptstyle{a}$};
\draw (1.35,-1.55) node {$\scriptscriptstyle{a}$};
\draw (1.6,-0.34) node {$\scriptscriptstyle{1}$};
\draw (3.7,-0.3) node {$\scriptscriptstyle{1}$};

\draw (5,-0.1) node {$=$};

\draw[thick] (6,-2)--(6.5,-1.5);
\draw[thick] (10,2)--(9.5,1.7);
\draw[thick] (9.5,-1.5)--(10,-2);
\draw (9.5,1.7)--(9.5,-1.5);

\draw[thick] (6,2)--(6.5,1.5);
\draw[thick] (6.5,1.5).. controls (8.7,0.7).. (9.5,-1.5);

\draw[thick] (7.7,0.5) .. controls (8.7,1.5)..(9.5,1.7);
\draw (7.8,-1)--(6.5,1.5);

\draw (6.5,-1.5)..controls (6.6,-0.2).. (7.7,0.5);
\draw [thick] (6.5,-1.5)..controls (7.4,-1.4).. (7.8,-1);
\draw [thick] (7.8,-1)..controls (8,-0.2).. (7.7,0.5);

\draw [snake=snake,segment amplitude=.2mm,segment length=2mm] (9,-2)--(9,2);

\filldraw[black] (6.6,-0.65) circle (1.8pt);

\draw (6.2,-2.2) node {$\scriptstyle{a+1}$};
\draw (5.9,2.2) node {$\scriptstyle{b}$};
\draw (8.9,-2.2) node {$\scriptstyle{1}$};
\draw (10.15,-2.2) node {$\scriptstyle{b}$};
\draw (10.2,2.2) node {$\scriptstyle{a+1}$};
\draw (6.75,0.15) node {$\scriptscriptstyle{1}$};
\draw (8.55,1.5) node {$\scriptscriptstyle{a}$};
\draw (6.35,-0.65) node {$\scriptscriptstyle{a}$};
\draw (7.35,-1.55) node {$\scriptscriptstyle{a}$};
\draw (7.6,-0.34) node {$\scriptscriptstyle{1}$};
\draw (9.7,-0.3) node {$\scriptscriptstyle{1}$};

\draw (11,-0.1) node {$=$};

\end{tikzpicture}
\end{center}

\begin{center}
\begin{tikzpicture}[scale=0.6]

\draw (5,-0.1) node {$=$};

\draw[thick] (6,-2)--(6.5,-1.5);
\draw[thick] (10,2)--(9.5,1.7);
\draw[thick] (9.5,-1.5)--(10,-2);
\draw (9.5,1.7)--(9.5,-1.5);

\draw[thick] (6,2)--(6.5,1.5);
\draw[thick] (6.5,1.5).. controls (8.7,0.7).. (9.5,-1.5);

\draw[thick] (7.7,0.5) .. controls (8.7,1.5)..(9.5,1.7);
\draw (6.55,-1)..controls (7.3,-0.3).. (6.5,1.5);

\draw (6.5,-1.5)..controls (6.6,-0.2).. (7.7,0.5);
\draw [thick] (6.5,-1.5)..controls (7.4,-1.4).. (7.8,-1);
\draw [thick] (7.8,-1)..controls (8,-0.2).. (7.7,0.5);

\draw [snake=snake,segment amplitude=.2mm,segment length=2mm] (9,-2)--(9,2);

\filldraw[black] (7.41,0.3) circle (1.8pt);

\draw (6.2,-2.2) node {$\scriptstyle{a+1}$};
\draw (5.9,2.2) node {$\scriptstyle{b}$};
\draw (8.9,-2.2) node {$\scriptstyle{1}$};
\draw (10.15,-2.2) node {$\scriptstyle{b}$};
\draw (10.2,2.2) node {$\scriptstyle{a+1}$};
\draw (6.5,-0.3) node {$\scriptscriptstyle{1}$};
\draw (8.55,1.5) node {$\scriptscriptstyle{a}$};
\draw (7.35,0.5) node {$\scriptscriptstyle{a}$};
\draw (7.55,-1.55) node {$\scriptscriptstyle{a-1}$};
\draw (7.27,-0.34) node {$\scriptscriptstyle{1}$};
\draw (9.7,-0.3) node {$\scriptscriptstyle{1}$};
\draw (6.35,-1.26) node {$\scriptscriptstyle{2}$};

\draw (13.3,-0.2) node {$\displaystyle{+\,\,\sum\limits_{i+j=a-1}}$};

\draw[thick] (16,-2)--(16.5,-1.5);
\draw[thick] (20,2)--(19.5,1.7);
\draw[thick] (19.5,-1.5)--(20,-2);
\draw (19.5,1.7)--(19.5,-1.5);

\draw[thick] (16,2)--(16.5,1.5);
\draw[thick] (16.5,1.5).. controls (18.7,0.7).. (19.5,-1.5);

\draw[thick] (17.7,0.5) .. controls (18.7,1.5)..(19.5,1.7);
\draw (16.5,-1.5)--(16.5,1.5);

\draw (17.8,-1)..controls (17.1,-0.2).. (17.7,0.5);
\draw [thick] (16.5,-1.5)..controls (17.4,-1.4).. (17.8,-1);
\draw [thick] (17.8,-1)..controls (18,-0.2).. (17.7,0.5);

\draw [snake=snake,segment amplitude=.2mm,segment length=2mm] (19,-2)--(19,2);

\filldraw[black] (16.5,-0.2) circle (1.8pt)
                        (17.27,-0.2) circle (1.8pt);

\draw (16.2,-2.2) node {$\scriptstyle{a+1}$};
\draw (15.9,2.2) node {$\scriptstyle{b}$};
\draw (18.9,-2.2) node {$\scriptstyle{1}$};
\draw (20.15,-2.2) node {$\scriptstyle{b}$};
\draw (20.2,2.2) node {$\scriptstyle{a+1}$};
\draw (16.32,-0.2) node {$\scriptscriptstyle{i}$};
\draw (18.55,1.5) node {$\scriptscriptstyle{a}$};
\draw (16.32,-0.95) node {$\scriptscriptstyle{1}$};
\draw (17.43,-0.77) node {$\scriptscriptstyle{1}$};
\draw (17.35,-1.55) node {$\scriptscriptstyle{a}$};
\draw (17.1,-0.2) node {$\scriptscriptstyle{j}$};
\draw (18.33,-0.5) node {$\scriptscriptstyle{a\!-\!1}$};
\draw (19.7,-0.3) node {$\scriptscriptstyle{1}$};

\end{tikzpicture}
\end{center}

\noindent The first equality above follows by applying the pitchfork lemma twice, while the second one follows from Dot Migration. In the first diagram of the last line we have also used one associativity of splitters, and that diagram is equal to zero since it contains a dotless digon. Finally, all summands in the last summation with $j<a-1$ are zero by Lemma \ref{vazlem} due to the presence of a digon with the edges of thicknesses $1$ and $a-1$, and so finally by (\ref{pomoc11}) the last diagram is equal to the second summand on the RHS of (\ref{a1b}), thus proving this case.\\

Finally, we are left with the general case. Now, we ''explode" the curly strand of thickness $c$ into $c$ thin lines. We shall "pass" each of these $c$ thin lines through the thick crossing by using (\ref{a1b}). There are two summands on the right-hand-side of (\ref{a1b}) - the first one which is just simple passing of a thin curly strand, and the second extra-term. In all of such obtained summands for $c$ thin strands, either all of them have simply passed through a
the thick crossing, or at least for some of the thin strands we have an extra term from (\ref{a1b}). In the latter case, let the first strand  that produces an extra term be the $j$-th one (counting from right to left). Then we have:

\begin{center}

\begin{tikzpicture}[scale=0.6]

\draw [thick] (0,-2)--(4,2);
\draw [thick] (4,-2)--(0,2);
\draw [thick,snake=snake,segment amplitude=.2mm,segment length=1mm] (2,-2)--(1,0);
\draw [thick,snake=snake,segment amplitude=.2mm,segment length=1mm] (2,2)--(1,0);

\draw (0.2,-2.1) node {$\scriptstyle{a}$};
\draw (2.2,-2.1) node {$\scriptstyle{c}$};
\draw (4.2,-2.1) node {$\scriptstyle{b}$};

\draw (5.5,-0.1) node {$=$};

\draw [thick] (7,-2)--(11,2);
\draw [thick] (11,-2)--(7,2);
\draw [thick,snake=snake,segment amplitude=.2mm,segment length=1mm] (9,-2)--(10,0);
\draw [thick,snake=snake,segment amplitude=.2mm,segment length=1mm] (9,2)--(10,0);

\draw (7.2,-2.1) node {$\scriptstyle{a}$};
\draw (9.2,-2.1) node {$\scriptstyle{c}$};
\draw (11.2,-2.1) node {$\scriptstyle{b}$};

\draw (12.85,-0.1) node {$\displaystyle{+\quad\sum\limits_{j=1}^c}$};

\draw [thick] (15,-2)--(19,2);
\draw [thick] (19,-2)--(15,2);
\draw (16.2,-0.8)--(16.2,0.8);
\draw (17.8,-0.8)--(17.8,0.8);
\draw [thick,snake=snake,segment amplitude=.2mm,segment length=.8mm] (17,-2)--(17.2,-1.5);
\draw [thick,snake=snake,segment amplitude=.2mm,segment length=.8mm] (17,2)--(17.2,1.5);

\draw [snake=snake,segment amplitude=.2mm,segment length=.8mm] (17.5,0)--(17.2,-1.5);
\draw [snake=snake,segment amplitude=.2mm,segment length=.8mm] (17.5,0)--(17.2,1.5);

\draw [thick,snake=snake,segment amplitude=.2mm,segment length=.8mm] (18.8,-0.8)--(17.2,-1.5);
\draw [thick,snake=snake,segment amplitude=.2mm,segment length=.8mm] (18.8,0.8)--(17.2,1.5);
\draw [thick,snake=snake,segment amplitude=.2mm,segment length=.8mm] (18.8,-0.8)--(18.8,0.8);

\draw [thick,snake=snake,segment amplitude=.2mm,segment length=.8mm] (15.5,-0.8)--(17.2,-1.5);
\draw [thick,snake=snake,segment amplitude=.2mm,segment length=.8mm] (15.5,0.8)--(17.2,1.5);
\draw [thick,snake=snake,segment amplitude=.2mm,segment length=.8mm] (15.5,-0.8)--(15.5,0.8);

\draw (15.2,-2.1) node {$\scriptstyle{a}$};
\draw (17.2,-2.1) node {$\scriptstyle{c}$};
\draw (17.2,2.1) node {$\scriptstyle{c}$};
\draw (19.2,-2.1) node {$\scriptstyle{b}$};
\draw (16.05,-0.1) node {$\scriptscriptstyle{1}$};
\draw (17.95,-0.1) node {$\scriptscriptstyle{1}$};

\draw (15.2,2.1) node {$\scriptstyle{b}$};
\draw (19.2,2.1) node {$\scriptstyle{a}$};

\draw (16.9,-0.75) node {$\scriptscriptstyle{j\!-\!1}$};
\draw (16.4,-1.45) node {$\scriptscriptstyle{\pi_j}$};

\draw (17.13,1.1) node {$\scriptscriptstyle{1}$};
\draw (15.05,-0.1) node {$\scriptscriptstyle{c-j}$};
\draw (19.35,0.1) node {$\scriptscriptstyle{j-1}$};

\filldraw[black] (16.5,-1.2) circle (2.5pt)
                        (17.32,-0.82) circle (2pt);

\end{tikzpicture}
\end{center}
Here we have denoted $\pi_j:={\pi_{_{ {K_{c-j,j}}   }}}$ and we have also used the associativity of splitters and (\ref{decdef}) to collect the different thin strands. Recall that $K_{p,q}$ denotes the partition $(q,q,\ldots,q)$ of length $p$. Now, we use the pitchfork lemma followed by the Thick R2 move on the curly strand of thickness $c-j$ and the left-most strand of thickness $1$, as well as on the curly strand of thickness $j-1$ and the right-most strand of thickness $1$. 

\begin{center}
\begin{tikzpicture}[scale=0.6]

\draw [thick] (0,-2)--(3,2);
\draw [thick] (3,-2)--(0,2);
\draw [thick,snake=snake,segment amplitude=.2mm,segment length=1mm] (1.5,-2)--(0.8,0);
\draw [thick,snake=snake,segment amplitude=.2mm,segment length=1mm] (1.5,2)--(0.8,0);

\draw (0.2,-2.1) node {$\scriptstyle{a}$};
\draw (1.7,-2.1) node {$\scriptstyle{c}$};
\draw (3.2,-2.1) node {$\scriptstyle{b}$};

\draw (4.1,-0.1) node {$=$};

\draw [thick] (5,-2)--(8,2);
\draw [thick] (8,-2)--(5,2);
\draw [thick,snake=snake,segment amplitude=.2mm,segment length=1mm] (6.5,-2)--(7.2,0);
\draw [thick,snake=snake,segment amplitude=.2mm,segment length=1mm] (6.5,2)--(7.2,0);

\draw (5.2,-2.1) node {$\scriptstyle{a}$};
\draw (6.7,-2.1) node {$\scriptstyle{c}$};
\draw (8.2,-2.1) node {$\scriptstyle{b}$};

\draw (11,-0.1) node {$\displaystyle{+\,\,\sum\limits_{j=1}^c\sum\limits_{k=0}^{c-j}\sum\limits_{r=0}^{j-1}}$};

\draw [thick] (14,-2)--(19,2);
\draw [thick] (19,-2)--(14,2);
\draw (14.45,-1.64)--(14.45,1.64);
\draw (18.25,-1.4)--(18.25,1.4);
\draw [thick,snake=snake,segment amplitude=.2mm,segment length=.8mm] (16.5,-2)--(16.5,-1.5);
\draw [thick,snake=snake,segment amplitude=.2mm,segment length=.8mm] (16.5,2)--(16.5,1.5);

\draw [snake=snake,segment amplitude=.2mm,segment length=.8mm] (17,0)--(16.5,-1.5);
\draw [snake=snake,segment amplitude=.2mm,segment length=.8mm] (17,0)--(16.5,1.5);

\draw [thick,snake=snake,segment amplitude=.2mm,segment length=.8mm] (17.8,-0.8)--(16.5,-1.5);
\draw [thick,snake=snake,segment amplitude=.2mm,segment length=.8mm] (17.8,0.8)--(16.5,1.5);
\draw [thick,snake=snake,segment amplitude=.2mm,segment length=.8mm] (17.8,-0.8)--(17.8,0.8);

\draw [thick,snake=snake,segment amplitude=.2mm,segment length=.8mm] (15,-0.45)--(16.5,-1.5);
\draw [thick,snake=snake,segment amplitude=.2mm,segment length=.8mm] (15,0.45)--(16.5,1.5);
\draw [thick,snake=snake,segment amplitude=.2mm,segment length=.8mm] (15,-0.45)--(15,0.45);

\draw (14.2,-2.1) node {$\scriptstyle{a}$};
\draw (16.7,-2.1) node {$\scriptstyle{c}$};
\draw (16.7,2.1) node {$\scriptstyle{c}$};
\draw (19.2,-2.1) node {$\scriptstyle{b}$};
\draw (14.25,-0.9) node {$\scriptscriptstyle{1}$};
\draw (18.4,-0.9) node {$\scriptscriptstyle{1}$};

\draw (14.2,2.1) node {$\scriptstyle{b}$};
\draw (19.2,2.1) node {$\scriptstyle{a}$};

\draw (16.4,-0.75) node {$\scriptscriptstyle{j\!-\!1}$};
\draw (16.1,-1.55) node {$\scriptscriptstyle{{}_{\pi_{j}}}$};
\draw (15.77,0) node {$\scriptscriptstyle{{}_{\varepsilon_{{}_{c\!-\!j\!-\!k}}}}$};
\draw (17.48,0) node {$\scriptscriptstyle{{}_{\varepsilon_{r}}}$};
\draw (14.25,0) node {$\scriptscriptstyle{k}$};
\draw (18.85,0) node {$\scriptscriptstyle{j\!-\!1\!-\!r}$};

\draw (16.55,0.9) node {$\scriptscriptstyle{1}$};
\draw (15.6,1.27) node {$\scriptscriptstyle{c\!-\!j}$};
\draw (17.5,1.27) node {$\scriptscriptstyle{j\!-\!1}$};

\filldraw[black] (16.2,-1.3) circle (2pt)
                        (14.45,0) circle (1.8pt)                        
                        (18.25,0) circle (1.8pt)
                        (15,0) circle (2.5pt)
                        (17.8,0) circle (2.5pt)
                        (16.72,-0.82) circle (2pt);

\end{tikzpicture}
\end{center}

Now, as in the proof of the Thick R2 move, the last diagram can be nonzero only for $r=0$, since otherwise we would have two thin curly strands in the "exploded" thick one, both with $j-1$ dots and therefore equal to zero by antisymmetry. Therefore we have:

\begin{equation}\label{vaznaf}
\begin{tikzpicture}[scale=0.6]

\draw [thick] (0,-2)--(3,2);
\draw [thick] (3,-2)--(0,2);
\draw [thick,snake=snake,segment amplitude=.2mm,segment length=1mm] (1.5,-2)--(0.8,0);
\draw [thick,snake=snake,segment amplitude=.2mm,segment length=1mm] (1.5,2)--(0.8,0);

\draw (0.2,-2.1) node {$\scriptstyle{a}$};
\draw (1.7,-2.1) node {$\scriptstyle{c}$};
\draw (3.2,-2.1) node {$\scriptstyle{b}$};

\draw (4.1,-0.1) node {$=$};

\draw [thick] (5,-2)--(8,2);
\draw [thick] (8,-2)--(5,2);
\draw [thick,snake=snake,segment amplitude=.2mm,segment length=1mm] (6.5,-2)--(7.2,0);
\draw [thick,snake=snake,segment amplitude=.2mm,segment length=1mm] (6.5,2)--(7.2,0);

\draw (5.2,-2.1) node {$\scriptstyle{a}$};
\draw (6.7,-2.1) node {$\scriptstyle{c}$};
\draw (8.2,-2.1) node {$\scriptstyle{b}$};

\draw (11.2,-0.1) node {$\displaystyle{+\quad\sum\limits_{j=1}^c\sum\limits_{k=0}^{c-j}}$};

\draw [thick] (14,-2)--(19,2);
\draw [thick] (19,-2)--(14,2);
\draw (14.45,-1.64)--(14.45,1.64);
\draw (18.25,-1.4)--(18.25,1.4);
\draw [thick,snake=snake,segment amplitude=.2mm,segment length=.8mm] (16.5,-2)--(16.5,-1.5);
\draw [thick,snake=snake,segment amplitude=.2mm,segment length=.8mm] (16.5,2)--(16.5,1.5);

\draw [thick,snake=snake,segment amplitude=.2mm,segment length=.8mm] (17.8,-0.8)--(16.5,-1.5);
\draw [thick,snake=snake,segment amplitude=.2mm,segment length=.8mm] (17.8,0.8)--(16.5,1.5);
\draw [thick,snake=snake,segment amplitude=.2mm,segment length=.8mm] (17.8,-0.8)--(17.8,0.8);

\draw [thick,snake=snake,segment amplitude=.2mm,segment length=.8mm] (15,-0.45)--(16.5,-1.5);
\draw [thick,snake=snake,segment amplitude=.2mm,segment length=.8mm] (15,0.45)--(16.5,1.5);
\draw [thick,snake=snake,segment amplitude=.2mm,segment length=.8mm] (15,-0.45)--(15,0.45);

\draw (14.2,-2.1) node {$\scriptstyle{a}$};
\draw (16.7,-2.1) node {$\scriptstyle{c}$};
\draw (16.7,2.1) node {$\scriptstyle{c}$};
\draw (19.2,-2.1) node {$\scriptstyle{b}$};
\draw (14.25,-0.9) node {$\scriptscriptstyle{1}$};
\draw (18.4,-0.9) node {$\scriptscriptstyle{1}$};

\draw (14.2,2.1) node {$\scriptstyle{b}$};
\draw (19.2,2.1) node {$\scriptstyle{a}$};

\draw (16.1,-1.55) node {$\scriptscriptstyle{{}_{\pi_{j}}}$};
\draw (15.77,0) node {$\scriptscriptstyle{{}_{\varepsilon_{{}_{c\!-\!j\!-\!k}}}}$};
\draw (14.25,0) node {$\scriptscriptstyle{k}$};
\draw (18.75,0) node {$\scriptscriptstyle{j\!-\!1}$};

\draw (15.6,1.27) node {$\scriptscriptstyle{c\!-\!j}$};
\draw (17.35,1.3) node {$\scriptscriptstyle{j}$};

\filldraw[black] (16.2,-1.3) circle (2pt)
                        (14.45,0) circle (1.8pt)                        
                        (18.25,0) circle (1.8pt)
                        (15,0) circle (2.5pt);

\end{tikzpicture}
 \end{equation}

We can now iterate this formula, by passing the thick curly strand of thickness $c\!-\!j$ through the thick crossing by using analogous formula (\ref{vaznaf}). In such a way we get (we are also applying dot slides):

\begin{center}
\begin{tikzpicture}[scale=0.6]

\draw [thick] (0,-1.5)--(3,1.5);
\draw [thick] (3,-1.5)--(0,1.5);
\draw [thick,snake=snake,segment amplitude=.2mm,segment length=1mm] (1.5,-1.5)--(0.8,0);
\draw [thick,snake=snake,segment amplitude=.2mm,segment length=1mm] (1.5,1.5)--(0.8,0);

\draw (0.2,-1.6) node {$\scriptstyle{a}$};
\draw (1.7,-1.6) node {$\scriptstyle{c}$};
\draw (3.2,-1.6) node {$\scriptstyle{b}$};

\draw (8,-0.1) node {$\displaystyle{=\quad\sum\limits_{x=0}^c \,\,\sum\limits_{1\le j_1<\cdots<j_x\le c} \,\,\sum\limits_{
k_1,\ldots,k_x} }$};

\draw (17,0) node {\quad};





\end{tikzpicture}
 \end{center}

\begin{equation}\label{velikaf}
\begin{tikzpicture}[scale=0.7]

\draw (5,0) node {$\quad$};
\draw[thick] (10,0.7778)--(20,3);
\draw[thick] (10,3)--(19,0)--(20,-0.333);

\draw (11,1)--(11,2.666);
\draw (10.25,0.823)--(10.25,2.9);
\draw (12.35,1.3)--(12.35,2.22);
\draw (13.25,1.5)--(13.25,1.94);

\draw (11.75,1.85) node {$\cdots$};

\draw [snake=snake,segment amplitude=.2mm,segment length=1mm] (15,-3.2)--(15,3.4);

\draw [snake=snake,segment amplitude=.2mm,segment length=1mm] (15,-2.8)--(19.5,-1.5);
\draw [snake=snake,segment amplitude=.2mm,segment length=1mm] (15,-2)--(17.6,-1.2);
\draw [snake=snake,segment amplitude=.2mm,segment length=1mm] (15,-0.5)--(16,0.3);
\draw [snake=snake,segment amplitude=.2mm,segment length=1mm] (15,0.4)--(15.5,0.75);

\draw [snake=snake,segment amplitude=.2mm,segment length=1mm] (19.5,3)--(19.5,-1.5);
\draw [snake=snake,segment amplitude=.2mm,segment length=1mm] (17.6,2.8)--(17.6,-1.2);
\draw [snake=snake,segment amplitude=.2mm,segment length=1mm] (16,2.6)--(16,0.3);
\draw [snake=snake,segment amplitude=.2mm,segment length=1mm] (15.5,2.5)--(15.5,0.75);

\draw [snake=snake,segment amplitude=.2mm,segment length=1mm] (19.5,3)--(15,3.2);
\draw [snake=snake,segment amplitude=.2mm,segment length=1mm] (17.6,2.8)--(15,3.1);
\draw [snake=snake,segment amplitude=.2mm,segment length=1mm] (16,2.6)--(15,2.9);
\draw [snake=snake,segment amplitude=.2mm,segment length=1mm] (15.5,2.5)--(15,2.8);

\draw (15.7,1.06)--(15.7,2.01);
\draw (16.4,0.85)--(16.4,2.2);
\draw (19.1,-0.033)--(19.1,2.78);
\draw (19.75,-0.25)--(19.75,2.95);

\draw (18.51,1.52) node {\tiny{$\scriptscriptstyle{{j_2\!-\!j_{1}\!-\!1 }   }  $}};
\draw (20.2,1.6) node {$\scriptscriptstyle{j_1\!-\!1}$};


\draw (16.1,-1) node {$\cdots\cdots$};
\draw (17,1.5) node {${\cdots}$};

\filldraw[black] (15,-2.6) circle (2pt)
               (15,-2.2) circle (2pt)
(15,-1.8) circle (2pt)
(15,-1.4) circle (2pt)
(15,-0.3) circle (2pt)
(15,0.1) circle (2pt)
(15,0.6) circle (2pt)
(15,1) circle (2pt);

\draw (14.18,-2.2) node {\tiny{$\scriptscriptstyle{\varepsilon_{{}_{c\!-\!j_1\!-\!k_1}}   }$}};
\draw (14.18,-1.3) node {\tiny{$\scriptscriptstyle{\varepsilon_{{}_{c\!-\!j_2\!-\!k_2}}   }$}};
\draw (13.95,0.1) node {\tiny{$\scriptscriptstyle{\varepsilon_{{}_{c\!-\!j_{x\!-\!1}\!-\!k_{x\!-\!1}}}   }$}};
\draw (14.18,1.1) node {\tiny{$\scriptscriptstyle{\varepsilon_{{}_{c\!-\!j_x\!-\!k_x}}   }$}};

\draw (14.18,-2.7) node {\tiny{$\scriptscriptstyle{\pi_{{}_{K_{c\!-\!j_1,j_1}}}   }$}};
\draw (13.92,-1.7) node {\tiny{$\scriptscriptstyle{\pi_{{}_{K_{c\!-\!j_2,j_2\!-\!j_1}}}   }$}};
\draw (13.4,-0.25) node {\tiny{$\scriptscriptstyle{\pi_{{}_{K_{c\!-\!j_{x\!-\!1},j_{x\!-\!1}\!-\!j_{x\!-\!2} }}}   }$}};
\draw (13.55,0.67) node {\tiny{$\scriptscriptstyle{\pi_{{}_{K_{c\!-\!j_{x\!-\!1},j_{x}\!-\!j_{x\!-\!1} }}}   }$}};

\filldraw[black] (11,1.7) circle (1.4pt)
                  (10.25,1.7) circle (1.4pt)
                    (12.35,1.7) circle (1.4pt)
                 (13.25,1.7) circle (1.4pt);

\draw (10.07,1.6) node {{\tiny{$\scriptscriptstyle{k_1}$}}};
\draw (10.8,1.6) node {{\tiny{$\scriptscriptstyle{k_2}$}}};
\draw (12.01,1.6) node {{\tiny{$\scriptscriptstyle{k_{x\!-\!1}}$}}};
\draw (13,1.7) node {{\tiny{$\scriptscriptstyle{k_x}$}}};

\filldraw[black] (15.7,1.7) circle (1.4pt)
                       (16.4,1.7) circle (1.4pt)
                    (19.75,1.7) circle (1.4pt)
                       (19.1,1.7) circle (1.4pt);

\draw[very thin,->] (17.5,3.8)--(15.8,1.75);
\draw[very thin,->] (18.5,3.4)--(16.5,1.75);

\draw (17.5,4) node {\tiny{$\scriptscriptstyle{j_x\!-\!j_{x\!-\!1}\!-\!1}$}};
\draw (18.77,3.55) node {\tiny{$\scriptscriptstyle{j_{x\!-\!1}\!-\!j_{x\!-\!2}\!-\!1}$}};

\draw (18,-2.2) node {$\scriptscriptstyle{j_1}$};
\draw (17,-1.8) node {\tiny{$\scriptscriptstyle{j_2\!-\!j_1}$}};
\draw (16.38,-0.1) node {\tiny{$\scriptscriptstyle{j_{x\!-\!1}\!-\!j_{x\!-\!2}}$}};
\draw (14.55,2.3) node {\tiny{$\scriptscriptstyle{c\!-\!j_x}$}};

\draw (10,0.5) node {$\scriptstyle{a}$};
\draw (20.1,2.9) node {$\scriptstyle{a}$};
\draw (10.2,3.2) node {$\scriptstyle{b}$};
\draw (15.1,-3.4) node {$\scriptstyle{c}$};
\draw (15.1,3.5) node {$\scriptstyle{c}$};
\draw (20,-0.5) node {$\scriptstyle{b}$};

\end{tikzpicture}
 \end{equation}

First of all, we note that we can assume that the summation is over all integers $k_i$, since the corresponding elementary symmetric polynomial $\varepsilon_{c-j_i-k_i}$ is equal to zero when $j_i+k_i$ is larger than $c$.
 
Now we apply the thick R2 move (after two pitchfork lemmas) on the curly strand of thickness $j_1$ and thin strand with $j_2-j_1-1$ dots.

\begin{equation}\label{vaf1}
\begin{tikzpicture}[scale=0.7]
\draw (5.5,1.6) node {$\displaystyle{\quad\sum\limits_{x=0}^c \,\,\sum\limits_{1\le j_1<\cdots<j_x\le c} \,\,\sum\limits_{
k_1,\ldots,k_x} \,\,\sum\limits_{r+s=j_1} }$};

\draw[thick] (10,0.7778)--(20,3);
\draw[thick] (10,3)--(19,0)--(20,-0.333);

\draw (11,1)--(11,2.666);
\draw (10.25,0.823)--(10.25,2.9);
\draw (12.35,1.3)--(12.35,2.22);
\draw (13.25,1.5)--(13.25,1.94);

\draw (11.75,1.85) node {$\cdots$};

\draw [snake=snake,segment amplitude=.2mm,segment length=1mm] (15,-3.2)--(15,3.4);

\draw [snake=snake,segment amplitude=.2mm,segment length=1mm] (15,-2.8)--(18.7,-1.7);
\draw [snake=snake,segment amplitude=.2mm,segment length=1mm] (15,-2)--(17.6,-1.2);
\draw [snake=snake,segment amplitude=.2mm,segment length=1mm] (15,-0.5)--(16,0.3);
\draw [snake=snake,segment amplitude=.2mm,segment length=1mm] (15,0.4)--(15.5,0.75);

\draw [snake=snake,segment amplitude=.2mm,segment length=1mm] (18.7,3.1)--(18.7,-1.7);
\draw [snake=snake,segment amplitude=.2mm,segment length=1mm] (17.6,2.8)--(17.6,-1.2);
\draw [snake=snake,segment amplitude=.2mm,segment length=1mm] (16,2.6)--(16,0.3);
\draw [snake=snake,segment amplitude=.2mm,segment length=1mm] (15.5,2.5)--(15.5,0.75);

\draw [snake=snake,segment amplitude=.2mm,segment length=1mm] (18.7,3.1)--(15,3.2);
\draw [snake=snake,segment amplitude=.2mm,segment length=1mm] (17.6,2.8)--(15,3.1);
\draw [snake=snake,segment amplitude=.2mm,segment length=1mm] (16,2.6)--(15,2.9);
\draw [snake=snake,segment amplitude=.2mm,segment length=1mm] (15.5,2.5)--(15,2.8);

\draw (15.7,1.06)--(15.7,2.01);
\draw (16.4,0.85)--(16.4,2.2);
\draw (19.3,-0.133)--(19.3,2.88);
\draw (19.75,-0.25)--(19.75,2.95);

\draw (20.71,3.25) node {\tiny{$\scriptscriptstyle{{j_2\!-\!j_{1}\!-\!1 }   }  $}};
\draw (20.2,1.7) node {$\scriptscriptstyle{j_1\!-\!1}$};


\draw (16.1,-1) node {$\cdots\cdots$};
\draw (17,1.5) node {${\cdots}$};

\filldraw[black] (15,-2.6) circle (2pt)
               (15,-2.2) circle (2pt)
(15,-1.8) circle (2pt)
(15,-1.4) circle (2pt)
(15,-0.3) circle (2pt)
(15,0.1) circle (2pt)
(15,0.6) circle (2pt)
(15,1) circle (2pt);

\draw (14.18,-2.2) node {\tiny{$\scriptscriptstyle{\varepsilon_{{}_{c\!-\!j_1\!-\!k_1}}   }$}};
\draw (14.18,-1.3) node {\tiny{$\scriptscriptstyle{\varepsilon_{{}_{c\!-\!j_2\!-\!k_2}}   }$}};
\draw (13.95,0.1) node {\tiny{$\scriptscriptstyle{\varepsilon_{{}_{c\!-\!j_{x\!-\!1}\!-\!k_{x\!-\!1}}}   }$}};
\draw (14.18,1.1) node {\tiny{$\scriptscriptstyle{\varepsilon_{{}_{c\!-\!j_x\!-\!k_x}}   }$}};

\draw (18.45,1.15) node {{$\scriptscriptstyle{\varepsilon_{{}_{r}}   }$}};
\draw (19.15,1.2) node {{$\scriptscriptstyle{s}   $}};

\draw (14.18,-2.7) node {\tiny{$\scriptscriptstyle{\pi_{{}_{K_{c\!-\!j_1,j_1}}}   }$}};
\draw (13.92,-1.7) node {\tiny{$\scriptscriptstyle{\pi_{{}_{K_{c\!-\!j_2,j_2\!-\!j_1}}}   }$}};
\draw (13.4,-0.25) node {\tiny{$\scriptscriptstyle{\pi_{{}_{K_{c\!-\!j_{x\!-\!1},j_{x\!-\!1}\!-\!j_{x\!-\!2} }}}   }$}};
\draw (13.55,0.67) node {\tiny{$\scriptscriptstyle{\pi_{{}_{K_{c\!-\!j_{x\!-\!1},j_{x}\!-\!j_{x\!-\!1} }}}   }$}};

\filldraw[black] (11,1.7) circle (1.4pt)
                  (10.25,1.7) circle (1.4pt)
                    (12.35,1.7) circle (1.4pt)
                 (13.25,1.7) circle (1.4pt);

\draw (10.07,1.6) node {{\tiny{$\scriptscriptstyle{k_1}$}}};
\draw (10.8,1.6) node {{\tiny{$\scriptscriptstyle{k_2}$}}};
\draw (12.01,1.6) node {{\tiny{$\scriptscriptstyle{k_{x\!-\!1}}$}}};
\draw (13,1.7) node {{\tiny{$\scriptscriptstyle{k_x}$}}};

\filldraw[black] (15.7,1.7) circle (1.4pt)
                       (16.4,1.7) circle (1.4pt)
                    (19.75,1.8) circle (1.4pt)
                    (19.3,1.3) circle (1.4pt)
                       (19.3,1.8) circle (1.4pt)
                       (18.7,1.26) circle (2pt);

\draw[very thin,->] (17.5,3.8)--(15.8,1.75);
\draw[very thin,->] (18.5,3.4)--(16.5,1.75);
\draw[very thin,->] (20.7,3.1)--(19.4,1.85);

\draw (17.5,4) node {\tiny{$\scriptscriptstyle{j_x\!-\!j_{x\!-\!1}\!-\!1}$}};
\draw (18.77,3.55) node {\tiny{$\scriptscriptstyle{j_{x\!-\!1}\!-\!j_{x\!-\!2}\!-\!1}$}};

\draw (18,-2.2) node {$\scriptscriptstyle{j_1}$};
\draw (17,-1.8) node {\tiny{$\scriptscriptstyle{j_2\!-\!j_1}$}};
\draw (16.38,-0.1) node {\tiny{$\scriptscriptstyle{j_{x\!-\!1}\!-\!j_{x\!-\!2}}$}};
\draw (14.55,2.3) node {\tiny{$\scriptscriptstyle{c\!-\!j_x}$}};

\draw (10,0.5) node {$\scriptstyle{a}$};
\draw (20.1,2.9) node {$\scriptstyle{a}$};
\draw (10.2,3.2) node {$\scriptstyle{b}$};
\draw (15.1,-3.4) node {$\scriptstyle{c}$};
\draw (15.1,3.5) node {$\scriptstyle{c}$};
\draw (20,-0.5) node {$\scriptstyle{b}$};

\end{tikzpicture}
\end{equation}

We will show that from all the summands appearing in this thick R2 move only the one with no dots on the thick curly line of thickness $j_1$ can be nonzero. Indeed, otherwise we move the thick dots through splitters by using \cite[formulas (2.68) and (2.63)]{thick}:

\begin{center}
\begin{tikzpicture}[scale=0.5]

\draw [thick] (0,-1.5)--(0,0);
\draw [thick] (0,0)--(-1,1.5);
\draw [thick] (0,0)--(1,1.5);

\draw (0.4,-1.7) node {$\scriptstyle{a+b}$};
\draw (-0.8,1.6) node {$\scriptstyle{a}$};
\draw (1.15,1.6) node {$\scriptstyle{b}$};

\filldraw[black] (0,-0.75) circle (2.5pt);

\draw (3.3,-0.25) node {$=\,\sum\limits_{\alpha,\beta} \,\, c_{\alpha,\beta}^{\gamma}$};

\draw [thick] (7,-1.5)--(7,0);
\draw [thick] (7,0)--(6,1.5);
\draw [thick] (7,0)--(8,1.5);

\draw (7.4,-1.7) node {$\scriptstyle{a+b}$};
\draw (6.2,1.6) node {$\scriptstyle{a}$};
\draw (8.15,1.6) node {$\scriptstyle{b}$};

\filldraw[black] (7.6,0.9) circle (2.5pt);
\filldraw[black] (6.4,0.9) circle (2.5pt);

\draw (-0.4,-0.88) node {$\scriptscriptstyle{\pi_{{}_{\gamma}}}$};
\draw (8.05,0.65) node {$\scriptscriptstyle{\pi_{{}_{\beta}}}$};
\draw (6.3,0.58) node {$\scriptscriptstyle{\pi_{{}_{\alpha}}}$};

\draw (11,0) node {and};

\draw [thick] (15,-1.5)--(15,0);
\draw [thick] (15,0)--(14,1.5);
\draw [thick] (15,0)--(16,1.5);

\draw (15.4,-1.7) node {$\scriptstyle{a+b}$};
\draw (14.2,1.6) node {$\scriptstyle{a}$};
\draw (16.15,1.6) node {$\scriptstyle{b}$};

\filldraw[black] (15,-0.75) circle (2.5pt);

\draw (17.7,-0.25) node {$=\,\sum\limits_{r+s=n} $};

\draw [thick] (21,-1.5)--(21,0);
\draw [thick] (21,0)--(20,1.5);
\draw [thick] (21,0)--(22,1.5);

\draw (21.4,-1.7) node {$\scriptstyle{a+b}$};
\draw (20.2,1.6) node {$\scriptstyle{a}$};
\draw (22.15,1.6) node {$\scriptstyle{b}$};

\filldraw[black] (21.6,0.9) circle (2.5pt);
\filldraw[black] (20.4,0.9) circle (2.5pt);

\draw (14.6,-0.8) node {$\scriptscriptstyle{\varepsilon_{_{n}}}$};
\draw (21.95,0.65) node {$\scriptscriptstyle{\varepsilon_{{}_{s}}}$};
\draw (20.15,0.69) node {$\scriptscriptstyle{\varepsilon_{{}_{r}}}$};

\end{tikzpicture}
\end{center}

Then, a part of the diagram containing thick curly edges of thicknesses $j_1$ and $j_2-j_1$, after moving up the thick dots labelled by $\varepsilon_{c-j_1-k_1}$ and $\pi_{{}_{K_{c-j_1,j_1}}}$ and dot sliding of $\varepsilon_r$, by associativity of splitters becomes:

\begin{center}
\begin{tikzpicture}[scale=0.65]

\draw [snake=snake,segment amplitude=.2mm,segment length=1mm] (0,-2)--(0,-1.8)--(0,-1)--(0,0);
\draw [snake=snake,segment amplitude=.2mm,segment length=1mm] (0,1.1)--(0,1.5)--(0,2);
\draw [snake=snake,segment amplitude=.2mm,segment length=1mm] (0,-1)--(1.5,-0.2)--(1.5,1.2)--(0,1.7);
\draw [snake=snake,segment amplitude=.2mm,segment length=1mm] (0,-1.8)--(3,-1)--(3,1.5)--(0,1.9);

\draw (0,0.7) node {$\vdots$};

\filldraw[black] (0,-1.6) circle (1.8pt)
                        (0,-1.2) circle (1.8pt)
                        (0,-0.7) circle (1.8pt)
                        (0,-0.3) circle (1.8pt)
                        (3,0) circle (1.8pt);

\draw (-1.05,-1.25) node {\tiny{$\scriptscriptstyle{\varepsilon_{{}_{c\!-\!j_1\!-\!k_1}}}$}};
\draw (-1.05,-0.2) node {\tiny{$\scriptscriptstyle{\varepsilon_{{}_{c\!-\!j_2\!-\!k_2}}}$}};

\draw (-1.15,-1.7) node {\tiny{$\scriptscriptstyle{\pi_{{}_{{K}_{c\!-\!j_1,j_1}}}}$}};
\draw (-1.35,-0.65) node {\tiny{$\scriptscriptstyle{\pi_{{}_{{K}_{c\!-\!j_2,j_2\!-\!j_1}}}}$}};

\draw (2.7,0) node {{$\scriptscriptstyle{\varepsilon_{{}_{r}}}$}};

\draw (0.1,-2.1) node {$\scriptstyle{c}$};
\draw (0.1,2.1) node {$\scriptstyle{c}$};
\draw (2,-1.6) node {\tiny{$\scriptscriptstyle{j_1}$}};
\draw (1.1,-0.77) node {\tiny{$\scriptscriptstyle{j_2\!-\!j_1}$}};

\draw (5,0) node {$\displaystyle{=\sum\limits_{y=0}^{c-j_1-k_1}}$};

\draw [snake=snake,segment amplitude=.2mm,segment length=1mm] (9.5,-2)--(9.5,-1.8)--(9.5,-1)--(9.5,0.5);
\draw [snake=snake,segment amplitude=.2mm,segment length=1mm] (9.5,1.6)--(9.5,2);
\draw [snake=snake,segment amplitude=.2mm,segment length=1mm] (9.5,-1.8)--(11.5,-1.5)--(11.5,-1.2);
\draw [snake=snake,segment amplitude=.2mm,segment length=1mm] (9.5,1.9)--(11.5,1.5)--(11.5,1.2);
\draw [snake=snake,segment amplitude=.2mm,segment length=1mm] (9.5,1.9)--(11.5,1.5)--(11.5,1.2);
\draw [snake=snake,segment amplitude=.2mm,segment length=1mm] (11.5,-1.2)--(10.5,-0.8)--(10.5,0.8)--(11.5,1.2);
\draw [snake=snake,segment amplitude=.2mm,segment length=1mm] (11.5,-1.2)--(13.5,-0.8)--(13.5,0.8)--(11.5,1.2);

\draw (9.5,1.2) node {$\vdots$};

\filldraw[black] (9.5,-1.5) circle (1.8pt)
                        (9.5,-1) circle (1.8pt)
                        (9.5,-0.5) circle (1.8pt)
                        (9.5,0) circle (1.8pt)
                        (10.5,-0.4) circle (1.8pt)
                        (10.5,0.38) circle (1.8pt)
                        (13.52,0) circle (1.8pt);

\draw (8.2,-1) node {\tiny{$\scriptscriptstyle{\varepsilon_{{}_{c\!-\!j_1\!-\!k_1\!-\!y}}}$}};
\draw (8.35,0.1) node {\tiny{$\scriptscriptstyle{\varepsilon_{{}_{c\!-\!j_2\!-\!k_2}}}$}};

\draw (8.3,-1.6) node {\tiny{$\scriptscriptstyle{\pi_{{}_{{K}_{c\!-\!j_2,j_1}}}}$}};
\draw (8.1,-0.4) node {\tiny{$\scriptscriptstyle{\pi_{{}_{{K}_{c\!-\!j_2,j_2\!-\!j_1}}}}$}};

\draw (10.9,0.37) node {{$\scriptscriptstyle{\varepsilon_{{}_{y}}}$}};
\draw (11.65,-0.45) node {\tiny{$\scriptscriptstyle{\pi_{{}_{{K}_{j_2\!-\!j_1,j_1}}}}$}};

\draw (13.85,0) node {{$\scriptscriptstyle{\varepsilon_{{}_{r}}}$}};

\draw (9.6,-2.1) node {$\scriptstyle{c}$};
\draw (9.6,2.1) node {$\scriptstyle{c}$};
\draw (10.6,-1.85) node {\tiny{$\scriptscriptstyle{j_2}$}};
\draw (12.5,1.2) node {\tiny{$\scriptscriptstyle{j_1}$}};
\draw (10.6,1.15) node {\tiny{$\scriptscriptstyle{j_2\!-\!j_1}$}};
\draw (10.6,1.9) node {\tiny{$\scriptscriptstyle{j_2}$}};

\end{tikzpicture}
\end{center}

Again, for $r>0$ due to antisymmetry the bubble can be nonzero only if $y=j_2-j_1$. But then the expression on the curly strands is symmetric in $k_1$ and $k_2$ (the only dependence on them is the product $\varepsilon_{{c-j_2-k_1}}\varepsilon_{{c-j_2-k_2}}$), whereas it is antisymmetric on the ordinary strands on the left of the diagram in (\ref{vaf1}), and therefore the whole diagram is zero. 

Hence, only the summand with $r=0$ survives in (\ref{vaf1}), and then by using Lemma \ref{pomocna}  the bubble appearing in the above picture becomes just the thick strand of thickness $j_2$ with the dot labelled by $\varepsilon_{{}_{y}}$. Finally, after moving down the dots decorated by $\varepsilon_{{y}}$ and $\varepsilon_{{c-j_1-k_1-y}}$ throughout the splitter, we get that (\ref{vaf1}) is equal to (note that $r=0$ and so $s=j_1$)

\begin{center}
\begin{tikzpicture}[scale=0.7]
\draw (6.5,1.6) node {$\displaystyle{\quad\sum\limits_{x=0}^c \,\,\sum\limits_{1\le j_1<\cdots<j_x\le c} \,\,\sum\limits_{
k_1,\ldots,k_x} \,\, }$};

\draw (5,0) node {$\quad$};
\draw[thick] (10,0.7778)--(20,3);
\draw[thick] (10,3)--(19,0)--(20,-0.333);

\draw (11,1)--(11,2.666);
\draw (10.25,0.823)--(10.25,2.9);
\draw (12.35,1.3)--(12.35,2.22);
\draw (13.25,1.5)--(13.25,1.94);

\draw (11.75,1.85) node {$\cdots$};

\draw [snake=snake,segment amplitude=.2mm,segment length=1mm] (15,-3.2)--(15,3.4);

\draw [snake=snake,segment amplitude=.2mm,segment length=1mm] (15,-2)--(18.1,-1.08);
\draw [snake=snake,segment amplitude=.2mm,segment length=1mm] (15,-0.5)--(16,0.3);
\draw [snake=snake,segment amplitude=.2mm,segment length=1mm] (15,0.4)--(15.5,0.75);

\draw [snake=snake,segment amplitude=.2mm,segment length=1mm] (18.1,2.75)--(18.1,-1.08);
\draw [snake=snake,segment amplitude=.2mm,segment length=1mm] (16,2.6)--(16,0.3);
\draw [snake=snake,segment amplitude=.2mm,segment length=1mm] (15.5,2.5)--(15.5,0.75);

\draw [snake=snake,segment amplitude=.2mm,segment length=1mm] (18.1,2.75)--(15,3.1);
\draw [snake=snake,segment amplitude=.2mm,segment length=1mm] (16,2.6)--(15,2.9);
\draw [snake=snake,segment amplitude=.2mm,segment length=1mm] (15.5,2.5)--(15,2.8);

\draw (15.7,1.09)--(15.7,2.05);
\draw (16.4,0.85)--(16.4,2.2);
\draw (19.1,-0.033)--(19.1,2.78);
\draw (19.75,-0.25)--(19.75,2.95);
\draw (17.7,0.43)--(17.7,2.48);

\draw (18.68,1.7) node {{$\scriptscriptstyle{{j_2-\!1 }   }  $}};
\draw (20.2,1.7) node {$\scriptscriptstyle{j_1\!-\!1}$};


\draw (16.1,-1) node {$\cdots\cdots$};
\draw (17,1.5) node {${\cdots}$};

\filldraw[black] (15,-2.6) circle (2pt)
               (15,-1.8) circle (2pt)
(15,-1.4) circle (2pt)
(15,-0.3) circle (2pt)
(15,0.1) circle (2pt)
(15,0.6) circle (2pt)
(15,1) circle (2pt);

\draw (14.18,-2.6) node {\tiny{$\scriptscriptstyle{\varepsilon_{{}_{c\!-\!j_1\!-\!k_1}}   }$}};
\draw (14.18,-1.3) node {\tiny{$\scriptscriptstyle{\varepsilon_{{}_{c\!-\!j_2\!-\!k_2}}   }$}};
\draw (13.95,0.1) node {\tiny{$\scriptscriptstyle{\varepsilon_{{}_{c\!-\!j_{x\!-\!1}\!-\!k_{x\!-\!1}}}   }$}};
\draw (14.18,1.1) node {\tiny{$\scriptscriptstyle{\varepsilon_{{}_{c\!-\!j_x\!-\!k_x}}   }$}};

\draw (14.12,-1.7) node {\tiny{$\scriptscriptstyle{\pi_{{}_{K_{c\!-\!j_2,j_2}}}   }$}};
\draw (13.4,-0.25) node {\tiny{$\scriptscriptstyle{\pi_{{}_{K_{c\!-\!j_{x\!-\!1},j_{x\!-\!1}\!-\!j_{x\!-\!2} }}}   }$}};
\draw (13.55,0.67) node {\tiny{$\scriptscriptstyle{\pi_{{}_{K_{c\!-\!j_{x\!-\!1},j_{x}\!-\!j_{x\!-\!1} }}}   }$}};

\filldraw[black] (11,1.7) circle (1.4pt)
                  (10.25,1.7) circle (1.4pt)
                    (12.35,1.7) circle (1.4pt)
                 (13.25,1.7) circle (1.4pt);

\draw (10.07,1.6) node {{\tiny{$\scriptscriptstyle{k_1}$}}};
\draw (10.8,1.6) node {{\tiny{$\scriptscriptstyle{k_2}$}}};
\draw (12.01,1.6) node {{\tiny{$\scriptscriptstyle{k_{x\!-\!1}}$}}};
\draw (13,1.7) node {{\tiny{$\scriptscriptstyle{k_x}$}}};

\filldraw[black] (15.7,1.7) circle (1.4pt)
                       (16.4,1.7) circle (1.4pt)
                    (19.75,1.8) circle (1.4pt)
                    (17.7,1.8) circle (1.4pt)
                  (19.1,1.8) circle (1.4pt);

\draw[very thin,->] (17.5,3.8)--(15.8,1.75);
\draw[very thin,->] (18.5,3.4)--(16.5,1.75);
\draw[very thin,->] (19.5,3.1)--(17.8,1.85);

\draw (17.5,4) node {\tiny{$\scriptscriptstyle{j_x\!-\!j_{x\!-\!1}\!-\!1}$}};
\draw (18.77,3.55) node {\tiny{$\scriptscriptstyle{j_{x\!-\!1}\!-\!j_{x\!-\!2}\!-\!1}$}};
\draw (19.5,3.2) node {\tiny{$\scriptscriptstyle{j_{3}\!-\!j_{2}\!-\!1}$}};

\draw (16.8,-1.75) node {\tiny{$\scriptscriptstyle{j_2}$}};
\draw (16.38,-0.1) node {\tiny{$\scriptscriptstyle{j_{x\!-\!1}\!-\!j_{x\!-\!2}}$}};
\draw (14.55,2.3) node {\tiny{$\scriptscriptstyle{c\!-\!j_x}$}};

\draw (10,0.5) node {$\scriptstyle{a}$};
\draw (20.1,2.9) node {$\scriptstyle{a}$};
\draw (10.2,3.2) node {$\scriptstyle{b}$};
\draw (15.1,-3.4) node {$\scriptstyle{c}$};
\draw (15.1,3.5) node {$\scriptstyle{c}$};
\draw (20,-0.5) node {$\scriptstyle{b}$};

\end{tikzpicture}
\end{center}

Now we can repeat the analogous thick R2 moves (preceded by the pitchfork lemma moves) in order to move all remaining thick curly strands to the left, starting with the strand of thickness $j_2$. In such a way we arrive at

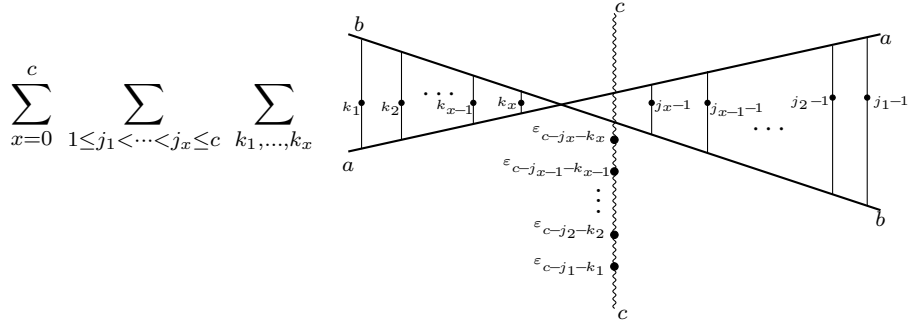
\begin{center}
\begin{tikzpicture}[scale=0.7]
\draw (6.5,1.6) node {$\displaystyle{\quad\sum\limits_{x=0}^c \,\,\sum\limits_{1\le j_1<\cdots<j_x\le c} \,\,\sum\limits_{
k_1,\ldots,k_x} \quad }$};

\draw (5,0) node {$\quad$};
\draw[thick] (10,0.7778)--(20,3);
\draw[thick] (10,3)--(19,0)--(20,-0.333);

\draw (11,1)--(11,2.666);
\draw (10.25,0.823)--(10.25,2.9);
\draw (12.35,1.3)--(12.35,2.22);
\draw (13.25,1.5)--(13.25,1.94);

\draw (11.75,1.85) node {$\cdots$};

\draw [snake=snake,segment amplitude=.2mm,segment length=1mm] (15,-2.2)--(15,3.4);

\draw (15.7,1.09)--(15.7,2.05);
\draw (16.75,0.75)--(16.75,2.3);
\draw (19.1,-0.033)--(19.1,2.78);
\draw (19.75,-0.25)--(19.75,2.95);

\draw (18.68,1.7) node {\tiny{$\scriptscriptstyle{{j_2-\!1 }   }  $}};
\draw (20.2,1.7) node {\tiny{$\scriptscriptstyle{j_1\!-\!1}$}};


\draw (17.9,1.2) node {${\cdots}$};

\filldraw[black] 
(15,-1.4) circle (2pt)
(15,-0.8) circle (2pt)
(15,0.4) circle (2pt)
(15,1) circle (2pt);

\draw (14.18,-1.4) node {\tiny{$\scriptscriptstyle{\varepsilon_{{}_{c\!-\!j_1\!-\!k_1}}   }$}};
\draw (14.18,-0.7) node {\tiny{$\scriptscriptstyle{\varepsilon_{{}_{c\!-\!j_2\!-\!k_2}}   }$}};
\draw (13.95,0.45) node {\tiny{$\scriptscriptstyle{\varepsilon_{{}_{c\!-\!j_{x\!-\!1}\!-\!k_{x\!-\!1}}}   }$}};
\draw (14.18,1.1) node {\tiny{$\scriptscriptstyle{\varepsilon_{{}_{c\!-\!j_x\!-\!k_x}}   }$}};

\draw (14.7,0) node {${\vdots}$};

\filldraw[black] (11,1.7) circle (1.4pt)
                  (10.25,1.7) circle (1.4pt)
                    (12.35,1.7) circle (1.4pt)
                 (13.25,1.7) circle (1.4pt);

\draw (10.07,1.6) node {{\tiny{$\scriptscriptstyle{k_1}$}}};
\draw (10.8,1.6) node {{\tiny{$\scriptscriptstyle{k_2}$}}};
\draw (12.01,1.6) node {{\tiny{$\scriptscriptstyle{k_{x\!-\!1}}$}}};
\draw (13,1.7) node {{\tiny{$\scriptscriptstyle{k_x}$}}};

\filldraw[black] (15.7,1.7) circle (1.4pt)
                       (16.75,1.7) circle (1.4pt)
                    (19.75,1.8) circle (1.4pt)
                  (19.1,1.8) circle (1.4pt);

\draw (16.1,1.65) node {\tiny{$\scriptscriptstyle{j_x\!-\!1}$}};
\draw (17.3,1.55) node {\tiny{$\scriptscriptstyle{j_{x\!-\!1}\!-\!1}$}};


\draw (10,0.5) node {$\scriptstyle{a}$};
\draw (20.1,2.9) node {$\scriptstyle{a}$};
\draw (10.2,3.2) node {$\scriptstyle{b}$};
\draw (15.15,-2.3) node {$\scriptstyle{c}$};
\draw (15.1,3.5) node {$\scriptstyle{c}$};
\draw (20,-0.5) node {$\scriptstyle{b}$};

\end{tikzpicture}
\end{center}

By using the associativity of splitters, we rewrite both groups of $x$ thin lines (the ones with labels $k_i$ and the other group with the labels $j_i-1$) as the exploded 
edges of thickness $x$. 
We do a similar thing as in the proof of the thick R2 move: first we note that due to the antisymmetry, the only nonzero summands are the ones where all $k_i$ are pairwise distinct, and we denote their decreasing ordering by $l_1>\ldots>l_x$. Let $\sigma\in S_x$ be a permutation such that $k_i=l_{\sigma_i}$, $i=1,\ldots,x$.
Then let $\alpha=(\alpha_1,\ldots,\alpha_x)$ and $\beta=(\beta_1,\ldots,\beta_x)$ be partitions defined by $\alpha_i+x-i=l_i$ and $\beta_{x+1-i}=j_i-i$, $i=1,\ldots,x$, respectively. 
Then the last expression becomes simply:

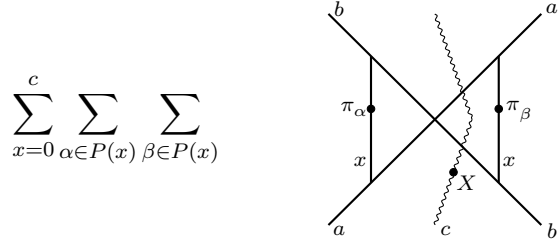
\begin{center}
\begin{tikzpicture}[scale=0.7]
\draw (10,0) node  {$\displaystyle{\sum_{x=0}^c\sum_{\alpha\in P(x)} \sum_{\beta\in P(x)}}$};

\draw [thick] (14,-2)--(18,2);
\draw [thick] (18,-2)--(14,2);
\draw [thick] (14.8,-1.2)--(14.8,1.2);
\draw [thick] (17.2,-1.2)--(17.2,1.2);
\draw [snake=snake,segment amplitude=.2mm,segment length=1mm] (16,-2)--(16.7,0);
\draw [snake=snake,segment amplitude=.2mm,segment length=1mm] (16,2)--(16.7,0);

\draw (14.2,-2.1) node {$\scriptstyle{a}$};
\draw (16.2,-2.1) node {$\scriptstyle{c}$};
\draw (18.2,-2.1) node {$\scriptstyle{b}$};
\draw (14.6,-0.8) node {$\scriptstyle{x}$};
\draw (17.4,-0.8) node {$\scriptstyle{x}$};
\draw (14.2,2.1) node {$\scriptstyle{b}$};
\draw (18.2,2.1) node {$\scriptstyle{a}$};

\filldraw[black] (14.8,0.2) circle (2pt)
                        (17.2,0.2) circle (2pt)
                        (16.35,-1) circle (2pt);
                        
\draw (14.5,0.2) node {$\scriptstyle{\pi_{{}_{\alpha}}}$};
\draw (17.6,0.15) node {$\scriptstyle{\pi_{{}_{\beta}}}$};
\draw (16.6,-1.2) node {$\scriptstyle{X}$};
                        
\end{tikzpicture}
\end{center}
where
\[
X= \sum_{\sigma\in S_x} \sgn(\sigma) \prod_{i=1}^x {\varepsilon_{c-\beta_{x+1-i}-i-\alpha_{\sigma_i}-x+\sigma_i}}.
\]
Finally, we have

\begin{eqnarray*}
&X=\sum\limits_{\sigma\in S_x} \sgn(\sigma)\prod\limits_{i=1}^x {\varepsilon_{{}_{c-\beta_{x+1-i}-i-\alpha_{{}_{\sigma_i}}-x+\sigma_i}}}= \\
&\quad = \det[\varepsilon_{{}_{c-\beta_{x+1-i}-i-\alpha_j-x+j}}]_{i,j=1}^x =\\
&\quad\quad=\det[\varepsilon_{{}_{(c-x-\beta_{x+1-i})-\alpha_j+j-i}}]_{i,j=1}^x =\\
&\quad\quad\quad\quad=\det [\varepsilon_{{}_{(K_{x}-\beta)_{i}-\alpha_j+j-i}}]_{i,j=1}^x=\pi_{\overline{K_{x}-\beta}{/}{\bar{\alpha}} }    .
\end{eqnarray*}
The last equality follows from the expression of the skew Schur polynomial as a determinant of the elementary symmetric polynomials \cite[formula (5.4), pp.70]{mac}. Finally, by definition, the last skew Schur polynomial can be written as follows:

\begin{eqnarray*}
&\pi_{\overline{K_{x}-\beta}/\bar{\alpha}}=\sum\limits_{\gamma} c_{\bar{\alpha},\bar{\gamma}}^{\overline{K_{x}-\beta}} \,\, \pi_{\bar{\gamma}}= \sum\limits_{\gamma} c_{{\alpha},{\gamma}}^{{K_{x}-\beta}} \,\, \pi_{\bar{\gamma}} =\\
&\quad\quad = \sum\limits_{\gamma} c_{{\alpha},{\gamma}}^{{K_{x}-\beta}} c_{\beta,K_{x}-\beta}^{K_{x}} \,\, \pi_{\bar{\gamma}} = \sum\limits_{\gamma} c_{{\alpha},\beta,{\gamma}}^{K_{x}}\,\, \pi_{\bar{\gamma}},
\end{eqnarray*}
\noindent which finishes our proof. \kraj

\section{Additional thick relations}\label{adtr}

The above thick R2 and R3 moves enable further relations for the general "passage" of a curly strand
through a thick edge, i.e. moving a strand labeled $c$ to the other side of the thick edge labeled $t$ in the picture below:

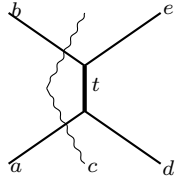
\begin{center}
 \begin{tikzpicture}[scale=0.5]

\draw [thick] (0,-2)--(2,-0.6);
\draw [thick] (2,0.6)--(4,2);
\draw [thick] (4,-2)--(2,-0.6);
\draw [thick] (2,0.6)--(0,2);
\draw [ultra thick] (2,0.6)--(2,-0.6);
\draw [snake=snake,segment amplitude=.2mm,segment length=2mm] (2,-2)--(1,0);
\draw [snake=snake,segment amplitude=.2mm,segment length=2mm] (2,2)--(1,0);

\draw (0.2,-2.1) node {$\scriptstyle{a}$};
\draw (2.2,-2.1) node {$\scriptstyle{c}$};
\draw (4.2,-2.1) node {$\scriptstyle{d}$};
\draw (0.2,2.1) node {$\scriptstyle{b}$};
\draw (4.2,2.1) node {$\scriptstyle{e}$};
\draw (2.3,0.1) node {$\scriptstyle{t}$};

\end{tikzpicture}
 \end{center}
 
with $a+b=d+e=t$.

Depending on the sign of the difference  $e-a=b-d$,  we have two slightly different formulas for this moving. Both
results are given in the following proposition.

\begin{proposition}
Let $a,b,c$ and $x$ be nonnegative integers. Then the following two relations hold in $\dot{\U}$:

\vskip 0.2cm
\begin{tikzpicture}[scale=0.5]

\draw [thick] (0,-2)--(2,-0.6);
\draw [thick] (2,0.6)--(4,2);
\draw [thick] (4,-2)--(2,-0.6);
\draw [thick] (2,0.6)--(0,2);
\draw [ultra thick] (2,0.6)--(2,-0.6);
\draw [snake=snake,segment amplitude=.2mm,segment length=2mm] (2,-2)--(1,0);
\draw [snake=snake,segment amplitude=.2mm,segment length=2mm] (2,2)--(1,0);

\draw (0.2,-2.1) node {$\scriptstyle{a}$};
\draw (2.2,-2.1) node {$\scriptstyle{c}$};
\draw (4.2,-2.1) node {$\scriptstyle{b+x}$};
\draw (0.2,2.1) node {$\scriptstyle{b}$};
\draw (4.2,2.1) node {$\scriptstyle{a+x}$};

\draw (5.5,-0.1) node {$=$};

\draw [snake=snake,segment amplitude=.2mm,segment length=2mm] (9,-2)--(10,0);
\draw [snake=snake,segment amplitude=.2mm,segment length=2mm] (9,2)--(10,0);


\draw [thick] (7,-2)--(11,2);
\draw [thick] (11,-2)--(7,2);
\draw (10.2,-1.2)--(10.2,1.2);

\draw (7.2,-2.1) node {$\scriptstyle{a}$};
\draw (9.2,-2.1) node {$\scriptstyle{c}$};
\draw (11.2,-2.1) node {$\scriptstyle{b+x}$};
\draw (10.4,-0.1) node {$\scriptscriptstyle{x}$};
\draw (7.2,2.1) node {$\scriptstyle{b}$};
\draw (11.2,2.1) node {$\scriptstyle{a+x}$};

\draw (12.5,-0.1) node {$+$};
\end{tikzpicture}

\begin{equation}\label{a1bn}
\begin{tikzpicture}[scale=0.4]
\draw (-2.5,0.3) node {$+$};
\draw (7.2,0.3) node {$\displaystyle{\sum\limits_{i=1}^{\min(a,b,c-x)}\sum\limits_{{\scriptsize{\begin{array}{cc}\alpha,\beta,\gamma\in\\
\in P(i,c\!-\!x\!-\!i)\end{array}}}}c_{\alpha,\beta,\gamma}^{\overbrace{{\scriptstyle (c-x-i,\ldots,c-x-i)}}^{i}} }$};
\end{tikzpicture}
\begin{tikzpicture}[scale=0.7]

\draw [thick] (14,-2)--(18,2);
\draw [thick] (18,-2)--(14,2);
\draw [thick] (14.8,-1.2)--(14.8,1.2);
\draw [thick] (17.2,-1.2)--(17.2,1.2);
\draw [snake=snake,segment amplitude=.2mm,segment length=1mm] (16,-2)--(16.7,0);
\draw [snake=snake,segment amplitude=.2mm,segment length=1mm] (16,2)--(16.7,0);

\draw (14.2,-2.1) node {$\scriptstyle{a}$};
\draw (16.2,-2.1) node {$\scriptstyle{c}$};
\draw (18.2,-2.1) node {$\scriptstyle{b+x}$};
\draw (14.6,-0.8) node {$\scriptstyle{i}$};
\draw (17.6,-0.8) node {$\scriptstyle{i+x}$};
\draw (14.2,2.1) node {$\scriptstyle{b}$};
\draw (18.2,2.1) node {$\scriptstyle{a+x}$};

\filldraw[black] (14.8,0.2) circle (2pt)
                        (17.2,0.2) circle (2pt)
                        (16.35,-1) circle (2pt);
                        
\draw (14.5,0.2) node {$\scriptstyle{\pi_{{}_{\alpha}}}$};
\draw (17.6,0.15) node {$\scriptstyle{\pi_{{}_{\beta}}}$};
\draw (16.6,-1.2) node {$\scriptstyle{\pi_{\bar{\gamma}}}$};
                        
\end{tikzpicture}
\end{equation}

\vskip 0.2cm

\begin{equation}\label{a1bnn}
\begin{tikzpicture}[scale=0.5]

\draw [thick] (0,-2)--(2,-0.6);
\draw [thick] (2,0.6)--(4,2);
\draw [thick] (4,-2)--(2,-0.6);
\draw [thick] (2,0.6)--(0,2);
\draw [ultra thick] (2,0.6)--(2,-0.6);
\draw [snake=snake,segment amplitude=.2mm,segment length=2mm] (2,-2)--(1,0);
\draw [snake=snake,segment amplitude=.2mm,segment length=2mm] (2,2)--(1,0);

\draw (0.2,-2.1) node {$\scriptstyle{a+x}$};
\draw (2.2,-2.1) node {$\scriptstyle{c}$};
\draw (4.2,-2.1) node {$\scriptstyle{b}$};
\draw (0.2,2.2) node {$\scriptstyle{b+x}$};
\draw (4.2,2.1) node {$\scriptstyle{a}$};

\draw (5.5,-0.1) node {$=$};

\draw (11.2,-0.15) node  {$\displaystyle{\sum_{i=0}^{\min(a,b,c)}\!\!\!\!\!\!\sum\limits_{{\scriptsize{\begin{array}{cc}\alpha,\beta,\gamma\in\\
\in P(i\!+\!x,c\!-\!i)\end{array}}}}\!\!\!\!\!\!c_{\alpha,\beta,\gamma}^{\overbrace{{\scriptstyle (c-i,\ldots,c-i)}}^{i+x}}}$};


\draw [thick] (16.6,-2)--(20.6,2);
\draw [thick] (20.6,-2)--(16.6,2);
\draw [thick] (17.4,-1.2)--(17.4,1.2);
\draw [thick] (19.8,-1.2)--(19.8,1.2);
\draw [snake=snake,segment amplitude=.2mm,segment length=1mm] (18.6,-2)--(19.3,0);
\draw [snake=snake,segment amplitude=.2mm,segment length=1mm] (18.6,2)--(19.3,0);

\draw (16.8,-2.1) node {$\scriptstyle{a+x}$};
\draw (18.8,-2.1) node {$\scriptstyle{c}$};
\draw (20.8,-2.1) node {$\scriptstyle{b}$};
\draw (16.8,-0.8) node {$\scriptstyle{i+x}$};
\draw (20,-0.8) node {$\scriptstyle{i}$};
\draw (16.8,2.2) node {$\scriptstyle{b+x}$};
\draw (20.8,2.1) node {$\scriptstyle{a}$};

\filldraw[black] (17.4,0.2) circle (2pt)
                        (19.8,0.2) circle (2pt)
                        (18.95,-1) circle (2pt);
                        
\draw (17,0.2) node {$\scriptstyle{\pi_{{}_{\alpha}}}$};
\draw (20.3,0.15) node {$\scriptstyle{\pi_{{}_{\beta}}}$};
\draw (19.25,-1.2) node {$\scriptstyle{\pi_{\bar{\gamma}}}$};
                        
\end{tikzpicture}
\end{equation}

\end{proposition}

\textbf{Proof:}

First we prove the formula (\ref{a1bn}). We note that by (\ref{OTE}):
\begin{equation}\label{up1}
\begin{tikzpicture}[scale=0.5]

\draw [thick] (0,-2)--(2,-0.6);
\draw [thick] (2,0.6)--(4,2);
\draw [thick] (4,-2)--(2,-0.6);
\draw [thick] (2,0.6)--(0,2);
\draw [ultra thick] (2,0.6)--(2,-0.6);
\draw [snake=snake,segment amplitude=.2mm,segment length=2mm] (2,-2)--(1,0);
\draw [snake=snake,segment amplitude=.2mm,segment length=2mm] (2,2)--(1,0);

\draw (0.2,-2.1) node {$\scriptstyle{a}$};
\draw (2.2,-2.1) node {$\scriptstyle{c}$};
\draw (4.2,-2.1) node {$\scriptstyle{b+x}$};
\draw (0.2,2.1) node {$\scriptstyle{b}$};
\draw (4.2,2.1) node {$\scriptstyle{a+x}$};

\draw (5.5,-0.1) node {$=$};

\draw [snake=snake,segment amplitude=.2mm,segment length=2mm] (9,-2)--(8,0);
\draw [snake=snake,segment amplitude=.2mm,segment length=2mm] (9,2)--(8,0);


\draw [thick] (7,-2)--(11,2);
\draw [thick] (11,-2)--(7,2);
\draw (10.2,-1.2)--(10.2,1.2);

\draw (7.2,-2.1) node {$\scriptstyle{a}$};
\draw (9.2,-2.1) node {$\scriptstyle{c}$};
\draw (11.2,-2.1) node {$\scriptstyle{b+x}$};
\draw (10.4,-0.1) node {$\scriptscriptstyle{x}$};
\draw (7.2,2.1) node {$\scriptstyle{b}$};
\draw (11.2,2.1) node {$\scriptstyle{a+x}$};

\end{tikzpicture}
\end{equation}

Now, by Thick R3 move (\ref{pr3}) we get that the RHS of (\ref{up1}) equals:

\begin{equation}\label{up2}
\begin{tikzpicture}[scale=0.5]

\draw (7,0) node {$\displaystyle{\sum_{i=0}^{\min(a,b,c)}\sum_{\alpha,\beta,\gamma\in P(i,c-i)}}c^{K_i}_{\alpha\beta\gamma}$};
\draw [thick](12,-3)--(13,-1.5);
\draw [thick](12,3)--(13,1.5);
\draw [thick](20,-3)--(17,-1.5);
\draw [thick](20,3)--(17,1.5);
\draw (18.5,2.25)--(18.5,-2.25);

\draw [thick](13,1.5)--(17,-1.5);
\draw [thick](13,-1.5)--(17,1.5);
\draw [thick](13,1.5)--(13,-1.5);
\draw [thick](17,1.5)--(17,-1.5);

\draw [snake=snake,segment amplitude=.2mm,segment length=2mm](15,-3)--(16,0);
\draw [snake=snake,segment amplitude=.2mm,segment length=2mm](16,0)--(15,3);

\draw (12.4,3) node {$\scriptstyle{b}$};
\draw (20.6,3) node {$\scriptstyle{a+x}$};
\draw (11.7,-3) node {$\scriptstyle{a}$};
\draw (20.6,-3) node {$\scriptstyle{b+x}$};
\draw (15.3,-3) node {$\scriptstyle{c}$};
\draw (18.75,1.2) node {$\scriptstyle{x}$};

\draw (12.8,1.2) node {$\scriptstyle{i}$};
\draw (17.2,1.2) node {$\scriptstyle{i}$};

\draw (14.2,1.2) node {$\scriptstyle{b-i}$};
\draw (14,-1.2) node {$\scriptstyle{a-i}$};

\draw (12.45,0) node {$\scriptstyle{\pi_{{}_{\alpha}}}$};
\draw (17.55,-0.05) node {$\scriptstyle{\pi_{{}_{\beta}}}$};
\draw (15.96,-2.1) node {$\scriptstyle{\pi_{{}_{\bar{\gamma}}}}$};

\filldraw[black] (13,0) circle (2pt) 
(17,0) circle (2pt)
(15.35,-2) circle (2pt);

\end{tikzpicture}
\end{equation}

\noindent where $K_i=(\overbrace{c-i,c-i,\ldots,c-i}^i)$, for $i>0$, and $K_0=0$. By  Associativity of splitters (Proposition \ref{asos}) and Lemma \ref{pomocna}, the formula (\ref{up2}) simplifies to:

\begin{equation}\label{up3}
\begin{tikzpicture}[scale=0.5]

\draw (7,0) node {$\displaystyle{\sum_{i=0}^{\min(a,b,c)}\sum_{\alpha,\beta,\gamma\in P(i,c-i)}}c^{K_i}_{\alpha\beta\gamma}$};
\draw [thick](12,-3)--(13,-1.5);
\draw [thick](12,3)--(13,1.5);
\draw [thick](18,-3)--(17,-1.5);
\draw [thick](18,3)--(17,1.5);

\draw [thick](13,1.5)--(17,-1.5);
\draw [thick](13,-1.5)--(17,1.5);
\draw [thick](13,1.5)--(13,-1.5);
\draw [thick](17,1.5)--(17,-1.5);

\draw [snake=snake,segment amplitude=.2mm,segment length=2mm](15,-3)--(16,0);
\draw [snake=snake,segment amplitude=.2mm,segment length=2mm](16,0)--(15,3);

\draw (12.4,3) node {$\scriptstyle{b}$};
\draw (18.6,3) node {$\scriptstyle{a+x}$};
\draw (11.7,-3) node {$\scriptstyle{a}$};
\draw (18.6,-3) node {$\scriptstyle{b+x}$};
\draw (15.3,-3) node {$\scriptstyle{c}$};

\draw (12.8,1.2) node {$\scriptstyle{i}$};
\draw (17.6,1.2) node {$\scriptstyle{i+x}$};

\draw (14.2,1.2) node {$\scriptstyle{b-i}$};
\draw (14,-1.2) node {$\scriptstyle{a-i}$};

\draw (12.45,0) node {$\scriptstyle{\pi_{{}_{\alpha}}}$};
\draw (18.15,-0.1) node {$\scriptstyle{\pi_{\beta/_{K_{i,x}}}}$};
\draw (15.96,-2.1) node {$\scriptstyle{\pi_{{}_{\bar{\gamma}}}}$};

\filldraw[black] (13,0) circle (2pt) 
(17,0) circle (2pt)
(15.35,-2) circle (2pt);

\end{tikzpicture}
\end{equation}

The last expression equals the RHS of (\ref{a1bn}), as wanted.\\

The formula (\ref{a1bnn}) is proved analogously: again, by Proposition \ref{OTE}  we have:\\

\begin{equation}\label{up4}
\begin{tikzpicture}[scale=0.5]

\draw [thick] (0,-2)--(2,-0.6);
\draw [thick] (2,0.6)--(4,2);
\draw [thick] (4,-2)--(2,-0.6);
\draw [thick] (2,0.6)--(0,2);
\draw [ultra thick] (2,0.6)--(2,-0.6);
\draw [snake=snake,segment amplitude=.2mm,segment length=2mm] (2,-2)--(1,0);
\draw [snake=snake,segment amplitude=.2mm,segment length=2mm] (2,2)--(1,0);

\draw (0.2,-2.1) node {$\scriptstyle{a+x}$};
\draw (2.2,-2.1) node {$\scriptstyle{c}$};
\draw (4.2,-2.1) node {$\scriptstyle{b}$};
\draw (0.2,2.1) node {$\scriptstyle{b+x}$};
\draw (4.2,2.1) node {$\scriptstyle{a}$};

\draw (5.5,-0.1) node {$=$};

\draw [snake=snake,segment amplitude=.2mm,segment length=2mm] (9,-2)--(7.4,0);
\draw [snake=snake,segment amplitude=.2mm,segment length=2mm] (9,2)--(7.4,0);


\draw [thick] (7,-2)--(11,2);
\draw [thick] (11,-2)--(7,2);
\draw (8.2,-0.8)--(8.2,0.8);

\draw (7.15,-2.2) node {$\scriptstyle{a+x}$};
\draw (9.2,-2.1) node {$\scriptstyle{c}$};
\draw (11.2,-2.1) node {$\scriptstyle{b}$};
\draw (8.4,-0.1) node {$\scriptscriptstyle{x}$};
\draw (7.15,2.2) node {$\scriptstyle{b+x}$};
\draw (11.2,2.1) node {$\scriptstyle{a}$};

\end{tikzpicture}
\end{equation}

By Pitchfork Lemma and Thick R2 move, the RHS of (\ref{up4}) becomes:

\begin{equation}\label{up6}
\begin{tikzpicture}[scale=0.5]

\draw (5.1,-0.1) node {${\displaystyle \sum\limits_{\varphi\in P(x,c)} }$};

\draw [snake=snake,segment amplitude=.2mm,segment length=2mm] (9,-2)--(8.2,0);
\draw [snake=snake,segment amplitude=.2mm,segment length=2mm] (9,2)--(8.2,0);


\draw [thick] (7,-2)--(11,2);
\draw [thick] (11,-2)--(7,2);
\draw (7.7,-1.3)--(7.7,1.3);

\draw (7.15,-2.2) node {$\scriptstyle{a+x}$};
\draw (9.05,2.2) node {$\scriptstyle{c}$};
\draw (11.2,-2.1) node {$\scriptstyle{b}$};
\draw (7.88,0.7) node {$\scriptscriptstyle{x}$};
\draw (7.15,2.2) node {$\scriptstyle{b+x}$};
\draw (11.2,2.1) node {$\scriptstyle{a}$};

\filldraw[black] (7.7,-0.2) circle (2pt) 
(8.85,-1.6) circle (2pt);

\draw (7.35,-0.2) node {$\scriptstyle{\pi_{\!{}_{\varphi}}}$};
\draw (9.26,-1.67) node {$\scriptstyle{\pi_{\!{}_{\widehat{\varphi}}}}$};

\end{tikzpicture}
\end{equation}

\noindent Recall, that by convention $\widehat{\varphi}=\overline{K_{x,c}-\varphi}$. Now, by Thick R3 move, the last expression is equal to

\begin{equation}\label{up8}
\begin{tikzpicture}[scale=0.5]

\draw (3.5,0) node {$\displaystyle{\sum_{\varphi\in P(x,c)}\sum_{i=0}^{\min(a,b,c)}\sum_{\alpha,\beta,\gamma\in P(i,c-i)}}c^{K_i}_{\alpha\beta\gamma}$};
\draw [thick](10,-3)--(13,-1.5);
\draw [thick](10,3)--(13,1.5);
\draw [thick](18,-3)--(17,-1.5);
\draw [thick](18,3)--(17,1.5);
\draw (11.5,2.25)--(11.5,-2.25);

\draw [thick](13,1.5)--(17,-1.5);
\draw [thick](13,-1.5)--(17,1.5);
\draw [thick](13,1.5)--(13,-1.5);
\draw [thick](17,1.5)--(17,-1.5);

\draw [snake=snake,segment amplitude=.2mm,segment length=2mm](15,-3)--(16,0);
\draw [snake=snake,segment amplitude=.2mm,segment length=2mm](16,0)--(15,3);

\draw (9.8,3.2) node {$\scriptstyle{b+x}$};
\draw (18.5,3) node {$\scriptstyle{a}$};
\draw (9.8,-3.15) node {$\scriptstyle{a+x}$};
\draw (18.5,-3) node {$\scriptstyle{b}$};
\draw (15.3,3) node {$\scriptstyle{c}$};
\draw (11.25,1.5) node {$\scriptstyle{x}$};

\draw (12.8,1.2) node {$\scriptstyle{i}$};
\draw (17.2,1.2) node {$\scriptstyle{i}$};

\draw (14.2,1.2) node {$\scriptstyle{b-i}$};
\draw (14,-1.2) node {$\scriptstyle{a-i}$};

\draw (13.5,0) node {$\scriptstyle{\pi_{{}_{\!\alpha}}}$};
\draw (11.1,0) node {$\scriptstyle{\pi_{{}_{\!\varphi}}}$};
\draw (15.7,-2.45) node {$\scriptstyle{\pi_{{}_{\!\widehat{\varphi}}}}$};
\draw (17.55,-0.05) node {$\scriptstyle{\pi_{{}_{\!\beta}}}$};
\draw (15.96,-1.6) node {$\scriptstyle{\pi_{{}_{\!\bar{\gamma}}}}$};

\filldraw[black] (13,0) circle (2pt) 
(11.5,0) circle (2pt) 
(17,0) circle (2pt)
(15.25,-2.25) circle (2pt)
(15.5,-1.5) circle (2pt);

\end{tikzpicture}
\end{equation}

Again by applying the Associativity of Splitters, and grouping the summand analogously as before, we obtain (\ref{a1bnn}).  

\kraj\\

\vskip 0.3cm

\textsc{CAMGSD,  Departamento de Matem\'atica, 
Instituto Superior T\'ecnico, 
Av. Rovisco Pais, 1049-001 Lisbon, Portugal,}

{and}

\textsc{Mathematical Institute SANU, Knez Mihailova 36, 11000 Beograd, Serbia.}

E-mail: {\tt{mstosic@math.ist.utl.pt}}

\end{document}